\documentclass[12pt,titlepage,reqno]{amsart}
\makeatletter
\newdimen\spfactor
\advance\topmargin by -1in
\advance\textwidth by 1.5in
\advance\oddsidemargin by -0.75in
\evensidemargin=\oddsidemargin 
\advance\textheight by \if@titlepage 1cm\else 1.8in\fi
\DeclareMathAlphabet\eusm{U}{eus}{m}{n}
\SetMathAlphabet\eusm{bold}{U}{eus}{b}{n}
\usepackage{amscd,amssymb,amsmath,cite}
\let\lbl\label
\def\subl#1{\subsection{}\label{#1}\gdef\lbl##1{\label{#1.##1}}\par}
\def\section{\def\@secnumfont{\bfseries}\@startsection{section}{1}%
  {\parindent}{.7\linespacing\@plus\linespacing}{.5\linespacing}%
  {\normalfont\bfseries\centering}}%
\def\sect{\def\@secnumfont{\bfseries}\@startsection{section}{1}%
  {\parindent}{.7\linespacing\@plus\linespacing}{.5\linespacing}%
  {\normalfont\bfseries}}%
\def\subsection{\def\@secnumfont{\bfseries}\@startsection{subsection}{2}%
  {\parindent}{.5\linespacing\@plus.7\linespacing}{-.5em}%
  {\normalfont\bfseries}}%
\def\appendix{\setcounter{equation}0\@ifnextchar*{\@app@star}{\@app@norm}}
\def\@app@star*#1{\def\thesection{A}\setcounter{subsection}0%
\sect*{Appendix. #1}}
\def\@app@norm#1{\if@firstapp\setcounter{section}{0}\@firstappfalse\fi%
\refstepcounter{section}\def\thesection{\Alph{section}}\sect*{%
Appendix~\thesection. #1}}
\let\ps@firstpage\ps@empty
\def\footnoterule{\kern-.4\p@
        \hrule\@width 12pc\kern11\p@\kern-.5\footnotesep}
\def\@settitle{\begin{center}%
\large\baselineskip20\p@\relax\bfseries\@title\end{center}}
\def\@setauthors{%
  \begingroup
  \trivlist
  \centering\footnotesize \@topsep30\p@\relax
  \advance\@topsep by -\baselineskip
  \item\relax
  \andify\authors
  \authors
  \endtrivlist
  \endgroup
}
\def\emptyset{\varnothing}

\def\phi{\varphi}

\def\End{\operatorname{End}}
\def\Hom{\operatorname{Hom}}
\def\Z{\mathbb Z}
\def\C{\mathbb C}
\def\R{\mathbb R}
\def\N{\mathbb N}
\def\Q{\mathbb Q}

\def\lie#1{\mathfrak{#1}}

\def\tensor{\otimes}

\def\union{\cup}

\def\intersection{\cap}
\def\underl{\mathaccent"7017}

\def\<#1>{\langle #1\rangle}
\def\Sl2{\lie{sl}(2)}

\def\iff{\ifmmode\;\Longleftrightarrow\;\else if and only if\ \fi}
\def\onlyif{\ifmmode\Longleftarrow\else only if\ \fi}
\let\@@@forall\forall

\def\forall{\ifmmode\,\@@@forall\, \else for all\fi}
\def\scal#1,#2 {\langle#1,#2^\vee\rangle}
\def\Scal#1,#2{(#1|#2)}

\newtheorem*{Prop}{Proposition}
\newtheorem*{Thm}{Theorem}
\newtheorem*{Lemma}{Lemma}
\newtheorem*{cor}{Corollary}

\theoremstyle{definition}
\newtheorem*{Rem}{Remark}
\newtheorem*{defi}{Definition}
\theoremstyle{remark}

\newtheorem*{example*}{Example}

\numberwithin{exercise}{section}
\numberwithin{example}{section}
\def\aref#1{\@ifundefined{r@#1}{\textup{(?)}}{\textup{(\@alph{\ref{#1}})}}}
\let\mak@titl@\maketitle
\def\maketitle{\baselineskip=0pt\lineskiplimit=0pt\normallineskiplimit=0pt%
\if@titlepage%
\global\advance\textheight by -1cm\global\advance\textheight by 1.8in\fi%
\mak@titl@\baselineskip=\spfactor\lineskiplimit=\spfactor%
\normallineskiplimit=\spfactor}

\newenvironment{enumer}{\begin{list}{{\hfill\rm(\roman{sbc})\hfill}}{%
\settowidth{\labelwidth}{(iii)}%
\usecounter{sbc}\leftmargin=\labelwidth\advance\leftmargin by\labelsep%
\rightmargin=\z@\itemindent=\z@}}
{\end{list}}
\newcounter{sbc}
\newenvironment{enumera}{\begin{list}{$\thesbc^\circ\mskip-6mu.$}%
{\settowidth{\labelwidth}{$5^\circ\mskip-6mu.$}\labelwidth=-\labelwidth
\advance\labelwidth by\labelsep%
\usecounter{sbc}\leftmargin=\z@\rightmargin=\z@\itemindent=\z@}}
{\end{list}}
\newenvironment{enumeraa}{\begin{list}{$\thesbc^\circ\mskip-6mu.$}{%
\settowidth{\labelwidth}{$5^\circ\mskip-6mu.$}%
\usecounter{sbc}\leftmargin=\labelwidth\advance\leftmargin by\labelsep%
\rightmargin=\z@\itemindent=\z@}}
{\end{list}}

{\endlist}
\newenvironment{enumspecm}[1]{\list{($#1\thesbc$)}{%
\settowidth{\labelwidth}{($#1{5}$)}%
\usecounter{sbc}\leftmargin=\labelwidth\advance\leftmargin by\labelsep%
\rightmargin=\z@\itemindent=\z@}}%
{\endlist}
\newenvironment{Proof}[1][\proofname]{\par
  \normalfont
  \topsep6pt plus6pt \trivlist
  \item[\hskip\labelsep\bfseries
    #1\ifnum\spacefactor>\@m\else{.}\fi]\ignorespaces
}{%
  \qed\endtrivlist
}
\@namedef{Proof*}{\Proof}
\@namedef{endProof*}{\endtrivlist}
\numberwithin{equation}{section}
\baselineskip=\spfactor
\lineskip=\spfactor
\normallineskip=\spfactor
\normallineskiplimit=\spfactor

\def\:{\,:\,}

\def\blambda{{\boldsymbol\lambda}}

\def\bmu{{\boldsymbol\mu}}

\def\gzero{\underl\lie g}
\def\hzero{\underl\lie h}

\def\to{\longrightarrow}
\def\mapsto{\longmapsto}

\def\qbinom#1#2{\genfrac[]{0pt}0{#1}{#2}}  
\makeatletter

\def\Sl2{\mathfrak{sl}_2}
\def\wt{\operatorname{wt}}
\def\sref#1#2#3{\ref{#1}~(${\mathrm #2}_{\ref{#2#3}}$)}
\def\ch{\operatorname{ch}}
\def\maj{\operatorname{desc}}
\def\tmaj{\widetilde{\maj}}
\def\phi{\varphi}
\def\epsilon{\varepsilon}
\def\Blm{B_\ell(m)}
\def\hatBlm#1{\widehat{\Blm}{}^{#1}}

\def\blambda{{\boldsymbol{\lambda}}}
\def\bmu{{\boldsymbol{\mu}}}
\def\ev{\operatorname{ev}}

\def\longtwoheadrightarrow{\DOTSB\protect\relbar\protect\joinrel%
\twoheadrightarrow}
\makeatother

\begin{document}
\title[]{\ \\[1cm]Littelmann's path crystal and combinatorics
of certain integrable~$\widehat{\lie{sl}_{\ell+1}}$~modules
of level zero}
\author[]{\normalsize by\\[12pt]{\bfseries Jacob Greenstein}\\[10pt]
Institut de Math\'ematiques de Jussieu (UMR 7586)\\
Universit\'e Pierre et Marie Curie\\
175 rue du Chevaleret, Plateau 7D\\
F-75013 Paris Cedex, France\\[10pt]
E-mail: greenste@math.jussieu.fr\\[1.5truein]
April 21, 2002\\
Revised: August 15, 2002
}
\thanks{\normalsize\hskip5mm\hyphenpenalty=10000 
This work was partially supported
by the Chateaubriand scholarship 
and by the EC TMR
network ``Algebraic Lie Representations'' grant no. ERB FMRX-CT97-0100}

\maketitle

\sect{Introduction}

Throughout the introduction, the base field is assumed to be~$\C$.

\subl{I5}
The aim of the present paper is to construct a subcrystal of
Littelmann's path crystal, whose formal character coincides with
that of a certain simple integrable module of level zero over 
the untwisted affine algebra associated with~$\lie{sl}_{\ell+1}$,
and to study the decomposition of the tensor product of that crystal
with a highest weight crystal.

Let~$\lie g$ be a symmetrizable Kac-Moody algebra with
a Cartan subalgebra~$\lie h$ and let~$\pi\subset\lie h^*$ 
be a set of simple roots of~$\lie g$. If~$\alpha\in\pi$,
denote the corresponding simple coroot by~$\alpha^\vee$
and let $x_{\pm \alpha}\in\lie g_{\pm\alpha}\setminus\{0\}$ be
the corresponding Chevalley generators of~$\lie g$.
Fix a weight lattice~$P(\pi)$ of~$\lie g$
and let~$P^+(\pi)=\{\lambda\in P(\pi)\:\alpha^\vee(\lambda)\ge0\}$
be the set of dominant weights. A $\lie g$~module~$M$ is called
integrable if~$M$ is a direct sum of its weight spaces~$M_\nu$, $\nu\in P(\pi)$
and the $x_{\pm \alpha}$ act locally nilpotently on~$M$ for all~$\alpha\in\pi$.
One can also define, in a similar way, a notion of
an integrable module over the quantized enveloping algebra~$U_q(\lie g)$
associated with~$\lie g$. If~$\dim M_\nu<\infty$ \forall~$\nu\in P(\pi)$,
we call~$M$ admissible and define its formal character by
$$
\operatorname{ch}M=\sum_{\nu\in P(\pi)}(\dim M_\nu)e^\nu.
$$

Given~$\lambda\in P^+(\pi)$, denote by~$V(\lambda)$ (respectively,
$V(-\lambda)$) the unique,
up to an isomorphism, highest (respectively, lowest) weight simple
integrable module over~$\lie g$ or over~$U_q(\lie g)$. 
Its character is given by the famous Weyl-Kac
formula (cf.~\cite[Chap.~10]{Kac92}) and, moreover, 
determines~$V(\lambda)$ up to an isomorphism. Another
important property of~$V(\lambda)$ is that it admits a crystal
basis and a canonical basis (cf.~\cite{JB,LB,Ka93}). 

\subl{I7}
Littelmann's path model provides a combinatorial realisation
of the crystal basis of~$V(\lambda)$, which reflects the
above properties of that module.
Namely,
let~$\mathbb P$ be the set of 
all piecewise-linear continuous paths~$b:[0,1]\to\R P(\pi)$ such that~$b(0)=0$ 
and~$b(1)\in P(\pi)$, where one identifies
$b$ and~$b'$ if~$b=b'$ up to a reparametrization.
After Littelmann~\cite{Li94,Li95}, one can 
endow~$\mathbb P$ with a structure
of a normal crystal, which will be henceforth referred to as
Littelmann's path crystal, by defining crystal operators~$e_\alpha,
f_\alpha$ \forall~$\alpha\in\pi$.  
Given a subcrystal~$B$
of~$\mathbb P$, define its formal character by
$$
\ch B=\sum_{b\in B} e^{b(1)}.
$$
Let~$\mathcal A$
be the associative monoid generated by the~$e_\alpha, f_\alpha\:\alpha\in\pi$. 
If~$\lambda\in P^+(\pi)$ 
and~$b_\lambda\in\mathbb P$ is a linear path connecting 
the origin with~$\lambda$, then the formal character of the 
subcrystal~$B(\lambda)=\mathcal Ab_\lambda$ of~$\mathbb P$ 
coincides with that of~$V(\lambda)$ ([\citen{Li94}]). Moreover,
$B(\lambda)$ provides a combinatorial model for the crystal basis 
of~$V(\lambda)$ and allows one to construct a standard monomial basis 
of~$V(\lambda)$~([\citen{Li98}]).

One of the fundamental properties of~$B(\lambda)$ is its independence
of the choice of~$b_\lambda$.
Namely, let~$b$ be a path in~$\mathbb P$, whose image
lies in the dominant Weyl chamber, that is~$\alpha^\vee(b(\tau))
\ge0$ for all~$\tau\in[0,1]$. If~$b'\in\mathbb P$ is another such a
path then, by the Isomorphism Theorem
of Littelmann ([\citen{Li95}, Theorem~7.1]), $\mathcal Ab$ is isomorphic
to~$\mathcal Ab'$ if and only if~$b(1)=b'(1)$. In particular, if
the image of~$b$ lies in the dominant Weyl chamber and~$b(1)=\lambda$,
then~$\mathcal Ab$ is isomorphic to~$B(\lambda)$. Thus,
similarly to~$V(\lambda)$, 
$B(\lambda)$ is uniquely determined, up to an isomorphism,
by its formal character.

\subl{I10}
If~$\lie g$ is not finite dimensional, there might exist simple admissible
integrable modules which are neither highest nor lowest weight.
The interest in this class of modules is due to
the observation that they occur as submodules in $\lie g$~modules 
$\Hom(V(\lambda),V(\mu))$, $\lambda,\mu\in P^+(\pi)$
(see for example~\cite[5.12]{J99} or~\cite[3.1]{JT}).
Namely, if~$V$ is a simple admissible integrable $\lie g$~module, 
denote by~$V^\#=\bigoplus_{\nu\in P(\pi)} V^*_\nu\subset V^*$
its graded dual. Then~$V^\#$ is also simple, admissible
and integrable and~$\Hom_{\lie g}(V^\#,\Hom(V(\lambda),V(\mu)))$
is isomorphic to the subspace $V^{\mu}_{\lambda-\mu}=\{
v\in V_{\lambda-\mu}\: x_\alpha^{\alpha^\vee(\mu)+1}v=0\}$.
In particular, if~$\lambda=\mu$, then~$V$ must have a non trivial
weight subspace of weight zero, which cannot occur in the highest or lowest
weight case. The embeddings of~$V^\#$ into~$\End V(\lambda)$ play a 
crucial role in the construction of KPRV determinants in the affine case
(cf.~\cite{J99,JT}). Thus one would like to be able to describe the
subspaces~$V^\mu_{0}$ or, more generally, $V^{\mu}_{\lambda-\mu}$. 
That problem is rather difficult for
modules, but is likely to
simplify significantly if one is able to pass to crystals.

\subl{I12}
Suppose now that $\lie g$ is a Kac-Moody algebra of an
untwisted affine type and denote by $\gzero$ its underlying finite 
dimensional simple
Lie algebra. By~\cite[Theorem~7.4]{Kac92},
$\lie g$ can be constructed
from $\gzero$ in the following way. Given a vector space~$V$, 
set~$L(V):=V\tensor\C[z,z^{-1}]$. Then~$\lie g$ is
the universal central extension of the semi-direct sum
of~$L(\gzero)$ and the one-dimensional space spanned by the Euler operator~$
D=z\frac{d}{dz}$. Let~$K$ be a central element of~$\lie g$. Then~$\lie g'=
[\lie g,\lie g]=L(\gzero)\oplus\C K$. A $\lie g$~module~$M$ is said
to be of level zero if~$K$ acts trivially on~$M$. One can easily see
that a highest or lowest weight module of level zero is necessarily
one-dimensional.

In the affine case, simple admissible integrable modules were classified
by V.~Chari and A.~Pressley ([\citen{C86,CP86,CP87}]). Moreover,
the modules of that type which are neither highest nor lowest weight
can be constructed as follows (cf.~\cite{CP86}).
For any~$\mathbf a:=(a_1,\dots,a_m)$, $a_i\in\C^\times$, define a
homomorphism of Lie algebras~$\ev_{\mathbf a}:\lie g'\to
\gzero^{\oplus m}$ by~$\ev_{\mathbf a}(K)=0$ 
and~$\ev_{\mathbf a}(x\tensor z^k)=
(a_1^k x,\dots,a_m^k x)$, \forall~$x\in\gzero$, $k\in\Z$. 
Let~$\mathbf V=(V_1,\dots,V_m)$ be a collection of finite-dimensional simple
$\gzero$~modules. Then~$V_1\tensor\cdots \tensor V_m$ 
is a simple~$\gzero^{\oplus m}$ module and we can endow it with a structure
of a~$\lie g'$~module by taking the pull-back by the 
homomorphism~$\ev_{\mathbf a}$. 
The resulting module is simple provided that all 
the~$a_j$ are distinct. 

Furthermore, the loop space~$L(V_1\tensor\cdots\tensor V_m)$
becomes a~$\lie g$~module, which we denote by $L(\mathbf V,\mathbf a)$,
if we set
$$
(x\tensor z^k)(v\tensor z^n)=(\ev_{\mathbf a}(x\tensor z^k))(v)
\tensor z^{k+n},\qquad D(v\tensor z^n)=nv\tensor z^n,
$$
\forall~$x\in\gzero$, $v\in V_1\tensor\cdots\tensor V_m$, $k,n\in\Z$.
If all the~$a_j$ are distinct, $L(\mathbf V,\mathbf a)$ is said to be
generic and
is completely
reducible. Simple submodules of modules of that
type exhausts all simple admissible integrable~$\lie g$~modules 
which are neither highest nor lowest weight modules.

Following~\cite[7.2]{J99}, 
we call these modules bounded 
since their weights satisfy the following condition.
By~\cite{Kac92}, there exists a non-degenerate symmetric bilinear form
$(\cdot,\cdot)$ on~$\lie h^*$ which is positive semidefinite on the
root lattice and may be assumed to be rational-valued on~$P(\pi)$.
A module~$M=\bigoplus_{\nu\in P(\pi)} M_\nu$ 
is called bounded if $(\nu,\nu)\le n$ for some~$n\in\N$ fixed and
\forall~$\nu\in P(\pi)$ such that~$M_\nu$ is non-trivial. 
If~$M$ is simple then the bound is actually attained
([\citen{J99}, 7.2])  
that is, there exists a weight~$\lambda\in P(\pi)$ of~$M$ called
maximal such that~$(\nu,\nu)\le (\lambda,\lambda)$ 
for all weights~$\nu$ of~$M$. For example, $V(\lambda)$ is always
bounded and~$\lambda$ is its maximal weight by~\cite[Proposition~11.4]{Kac92}.
One can show that a simple integrable module is admissible if
and only if it is bounded (cf.~\cite{J98}).

Formal characters of simple generic bounded modules 
were computed in~\cite{Gr99,Grxx}. It turns out that, unlike
the modules~$V(\lambda)$, these modules are not in general determined
by their formal characters up to an isomorphism. Besides, their
construction arises from the realisation of~$\lie g$ as a central
extension of a loop algebra, which is peculiar to 
Kac-Moody algebras of affine type. 
Thus, one should not expect that a combinatorial 
model similar to that of Littelmann for~$V(\lambda)$ exists
for an arbitrary simple admissible integrable module of level zero. 

\subl{I20}
Suppose now that~$a_i=\zeta^i$ where~$\zeta$ is an~$m$th primitive
root of unity and that~$V_1\cong \cdots\cong V_m\cong V$. Then
$L(\mathbf V,\mathbf a)$ becomes a direct sum of simple components
$L(V,m)^r$, $r=0,\dots,m-1$, where~$L(V,
m)^r$
is a cyclic submodule generated by~$v^{\tensor m}\tensor z^r$ and $v$ is
a highest weight vector of~$V$. The interest of this particular case
is due to the fact that the~$L(V,m)^r$ are determined by
their formal characters up to an isomorphism.

In the present paper we consider the case 
of~$\gzero\cong\lie{sl}_{\ell+1}$
and~$V$ isomorphic to the natural representation~$\C^{\ell+1}$ of~$\gzero$. 
Henceforth we will denote the corresponding modules~$L(V,m)^n$ 
by~$L(\ell,m;n)$. We show that $L(\ell,m;n)$ does admit a combinatorial
model in the framework of Littelmann's path crystal. Namely,
let~$\hzero$ be a Cartan subalgebra of~$\gzero$ and
$\omega$ be the highest weight of~$V$ with respect to~$\gzero$. 
Extend~$\omega$ to
the Cartan subalgebra~$\lie h=\hzero\oplus\C K\oplus\C D$ of~$\lie g$
by~$\omega(K)=\omega(D)=0$. 
Furthermore, let~$\delta\in\lie h^*$
be the unique element defined by the conditions~$\delta(D)=1$, 
$\delta|_{\hzero\oplus\C K}=0$. Then~$m\omega+n\delta$ is a maximal
weight of~$L(\ell,m;n)$. Needless to say, $m\omega+n\delta$ does not
lie in the dominant Weyl chamber.
Our main result is the following
\begin{Thm}
Let~$p_{\ell,m,n}$ be the linear path in~$\mathbb P$ connecting the
origin to~$m\omega+n\delta$. Then the formal character of the
subcrystal~$\hatBlm{n}=\mathcal Ap_{\ell,m,n}$ 
of~$\mathbb P$ 
equals the formal character of~$L(\ell,m;n)$.
\end{Thm}
As a byproduct, we obtain (Lemma~\ref{A20}) 
a nice combinatorial interpretation of
the dimensions of weight spaces of~$L(\ell,m;n)$. A similar result
holds for~$L(V^*,m)^n$, which is isomorphic to $L(\ell,m;m-n)^\#$. 

A natural question is how
the module~$L(\ell,m;n)$ is related to 
the crystal~$\hatBlm{n}$, apart from the
equality of their formal characters. It is shown
in~\cite{CG} that~$L(\ell,m;n)$ has a quantum analogue which, in turn,
admits a pseudo-crystal basis. The crystal $\hatBlm{n}$ provides
a combinatorial model for that basis (cf.~\cite[4.8--4.10]{CG}).

\subl{I30}
As a first application of the above result, 
we consider the decomposition of~$B(\lambda)\tensor \hatBlm{}$  
where~$\hatBlm{}=
\coprod_{n=0}^{m-1} \hatBlm{n}$ and the tensor product is understood
as concatenation of paths.
We obtain the following
\begin{Thm}
Let~$\lambda\in P(\pi)$ be a dominant weight which is not
a multiple of~$\delta$ and let
${\hatBlm{}}^\lambda$ be the set of paths~$b\in\hatBlm{}$ satisfying
$\alpha^\vee(\lambda+b(\tau))\ge0$ \forall~$\alpha\in\pi$ and
\forall~$\tau\in[0,1]$. 
Then the decomposition of the tensor product of~$B(\lambda)$ and~$\hatBlm{}$
is given by
$$
B(\lambda)\tensor \hatBlm{}\stackrel{\sim}{\to} \coprod_{b\in
{\hatBlm{}}^\lambda} B(\lambda+b(1)).
$$
\end{Thm}
We also obtain (Proposition~\ref{L60}) 
an explicit description of~${\hatBlm{}}^\lambda$ for the
case when~$\lambda$ is a fundamental weight. 

The above decomposition should be compared with the Decomposition rule
(cf.~\cite{Li94,Li95}), which generalizes the Littlewood-Richardson
rule, and with~\cite[Theorem~3.1]{HKKOTY}.
The main difference with the latter is that the crystal involved
in our situation is not finite. Besides, we consider an entirely
different framework, namely that of Littelmann's path crystal,
and our proof is not based on the theory of perfect crystals.
On the other hand, unlike that
of the Decomposition rule of~\cite{Li94,Li95}, the meaning of our decomposition
for modules is not yet understood. We expect, however, that it will
allow one to extract some information about embeddings of~$L(\ell,m;n)$
or its graded dual
into~$\Hom(V(\lambda),V(\mu))$, $\lambda,\mu\in P^+(\pi)$ discussed
in~\ref{I10}. 

\subsection*{Acknowledgements}
I am greatly indebted to A.~Joseph, whose lectures at
Universit\'e Pierre et Marie Curie (Paris VI) inspired
me to consider this problem in the first place, for acquainting me
with crystals, for his constant attention to my work and valuable advice.
An important part of the present paper was completed while I was 
visiting the MSRI, Berkeley, and
it is a great pleasure to thank the Institute and the organizers of
the programme ``Infinite-Dimensional Algebras and Mathematical Physics''
for their hospitality, the inspiring atmosphere and excellent working 
conditions.
I would like to express my gratitude to J.~Bernstein, 
V.~Chari, M.~Rosso, E.~Vasserot and S.~Zelikson for interesting
discussions and to R.~Rentschler for his hospitality and support. 
I am also obliged to the UMS MEDICIS, whose computational resources were 
used extensively in my work on the present paper. Finally, I would
like to thank the referee for numerous helpful suggestions.

\sect{Preliminaries}
In this section we recall 
the definition and some basic properties
of crystals and fix the notations
which will be used throughout the rest of the paper. Henceforth,
$\N$ stands for the set of non-negative integers and~$\N^+=\N\setminus\{0\}$.
The cardinality of a finite set~$S$ will be denoted by~$\# S$.

\subl{P10}
Let~$I$ be a finite index set and let~$A=(a_{ij})_{i,j\in I}$ be a
generalised Cartan matrix, that is, $a_{ii}=2$, $a_{ij}\in -\N$ if~$i\not=j$
and~$a_{ij}=0$ \iff~$a_{ji}=0$. We will assume that~$A$ is symmetrizable,
that is, there exist~$d_i\: i\in I$ such that the matrix~$(d_i a_{ij})_{i,j
\in I}$ is symmetric. 

Consider a triple~$(\lie h,\pi,\pi^\vee)$, where~$\lie h$ is
a~$\Q$-vector space, $\pi=\{\alpha_i\}_{i\in I}\subset\lie h^*$ 
and~$\pi^\vee=\{\alpha_i^\vee\}_{i\in I}$ is a linearly
independent subset of~$\lie h$. We would like to emphasize that $\pi$ is
{\em not} assumed to be linearly independent. We call
such a triple a realisation of~$A$ if~$a_{ij}=\alpha_i^\vee(\alpha_j)$, 
\forall~$i,j
\in I$. 
The realisation becomes unique,
up to an isomorphism, if we require both sets~$\pi$ and~$\pi^\vee$ to be
linearly independent and~$\dim\lie h=2 \#I-\operatorname{rk}A$ 
(cf.~\cite[Chap.~1]{Kac92}).

Given~$A$ and its realisation~$(\lie h,\pi,\pi^\vee)$, 
fix~$\Lambda_i\in\lie h^*$, $i\in I$ such
that~$\alpha_i^\vee(\Lambda_j)=\delta_{i,j}$ where~$\delta_{i,j}$ is
the Kronecker's symbol. Set~$P_0(\pi)=\bigoplus_{i\in I}\Z\Lambda_i$.
Complete the set $\{\Lambda_i\:i\in I\}$ to a basis of~$\lie h^*$,
and let~$P(\pi)$ be the free abelian group generated by that basis.

Endow~$\Z\union\{-\infty\}$ with a structure of an ordered semi-group
such that~$-\infty$ is the smallest element, $-\infty+n=-\infty$
\forall~$n\in\Z$ and~$\Z$ is given its natural order.
\begin{defi}[{cf.~\cite[Definition~1.2.1]{Ka93} or~\cite[5.2.1]{JB}}]
A crystal~$B$ is a set endowed with the maps~$e_i,f_i:B\to B\sqcup \{0\}$,
$\epsilon_i,\phi_i: B\to\Z\union\{-\infty\}$, $\wt: B\to P(\pi)$,
\forall~$i\in I$
which satisfy the following rules
\begin{enumspecm}{{\mathrm C}_}
\item\label{Cepphi} For any~$b\in B$, 
$\phi_i(b)=\epsilon_i(b)+\alpha_i^\vee(\wt b)$, 
\forall~$i\in I$.

\item\label{Cwteps} If~$b\in B$ and~$e_i b\in B$ (respectively,
$f_ib\in B$), then
$\wt e_i b=\wt b+\alpha_i$, $\epsilon_i(e_ib)=\epsilon_i(b)-1$
and~$\phi_i(e_ib)=\phi_i(b)+1$ (respectively, $\wt f_ib=
\wt b-\alpha_i$, $\epsilon_i(f_ib)=\epsilon_i(b)+1$
and~$\phi_i(f_ib)=\phi_i(b)-1$).

\item\label{Cinv} For $b,b'\in B$ and~$i\in I$,  $b' = e_i b$ \iff~$b=f_i b'$.

\item\label{Cinfty} If~$\phi_i(b)=-\infty$, then~$e_ib=f_ib=0$.
\end{enumspecm}
\end{defi}
\noindent
Given~$b\in B$, the value of~$\wt b$ is called the weight of~$b$. 

A crystal is said to be upper (respectively, lower) normal
if~$\epsilon_i(b)=\max\{ n\: e_i^n b\not=0\}$ (respectively,
$\phi_i(b)=\max\{ n\: f_i^n b\not=0\}$). A crystal is normal if it is
both upper and lower normal. 

\subl{P20}
Let~$B$ a crystal. For any~$\lambda\in P(\pi)$, set~$B_\lambda=\{b\in B\:
\wt b=\lambda\}$. If~$\# B_\lambda <\infty$ \forall~$\lambda\in P(\pi)$, 
one can define a formal character of~$B$ as
$$
\ch B=\sum_{b\in B}e^{\wt b}=\sum_{\lambda\in P(\pi)} \# B_\lambda e^\lambda.
$$
We say that~$\lambda\in P(\pi)$ is a weight of~$B$ if~$B_\lambda$ is non-empty.
Denote by~$\Omega(B)$ the set of all weights of~$B$.

\subl{P30}
Let~$B_1,\dots,B_n$ be crystals. The set~$B_1\times\cdots\times B_n$
can be endowed with a structure of a crystal which will be denoted
by~$B_1\tensor\cdots\tensor B_n$ and called the tensor product of crystals
$B_1,\dots, B_n$. The crystal maps are defined as follows 
(cf.~\cite[1.3]{Ka93}).

Given~$b=b_1\tensor\cdots \tensor b_n\in B_1\tensor\cdots\tensor B_n$,
define the Kashiwara functions~$b\mapsto r_k^i(b)\: i\in I$, $k\in\{1,
\dots, n\}$ by
$$
r_k^i(b)=\epsilon_i(b_k)-\sum_{1\le j<k} \alpha_i^\vee(\wt b_j).
$$
Then
\begin{enumspecm}{{\mathrm T}_}
\item\label{Teps} $\epsilon_i(b)=\max\{r^i_k(b)\: 1\le k\le n\}$.
\item\label{Twt} $\wt b=\sum_k \wt b_k$.
\item\label{Te} $e_i b=b_1\tensor \cdots\tensor  b_{r-1}\tensor e_i b_r
\tensor b_{r+1}\tensor \cdots \tensor b_n$, where~$r=\min\{k \:
r_k^i(b)=\epsilon_i(b)\}$, that is, $e_i$
acts in the leftmost place where the maximal value of~$r^i_k(b)$
is attained
\item\label{Tf} $f_i b=b_1\tensor \cdots\tensor b_{r-1}\tensor f_i b_r\tensor
b_{r+1}\tensor\cdots\tensor b_n$, where~$r=\max\{k\: r_k^i(b)=
\epsilon_i(b)\}$, that is, $f_i$ acts in
the rightmost place where the maximal value of~$r^i_k(b)$ is 
attained.
\end{enumspecm}
In the above we identify~$b_1\tensor\cdots\tensor b_{r-1}
\tensor 0\tensor b_{r+1}\tensor \cdots\tensor b_n$ with~$0$.
Since~$\phi_i(b)=\epsilon_i(b)+\alpha_i^\vee(\wt b)$, these rules
take a particularly nice form for the product~$B_1\tensor B_2$,
namely
\begin{align}
&e_i(b_1\tensor b_2)=\begin{cases}
e_i b_1\tensor b_2,&\text{if~$\phi_i(b_1)\ge \epsilon_i(b_2)$},\\
b_1\tensor e_ib_2,&\text{if~$\phi_i(b_1) < \epsilon_i(b_2)$}
\end{cases}
\lbl{10a}
\\
&f_i(b_1\tensor b_2)=\begin{cases}
f_i b_1\tensor b_2,&\text{if~$\phi_i(b_1)>\epsilon_i(b_2)$},\\
b_1\tensor f_ib_2,&\text{if~$\phi_i(b_1)\le \epsilon_i(b_2)$}
\end{cases}
\lbl{10b}
\end{align}
whilst~$\epsilon_i(b_1\tensor b_2)=\max\{\epsilon_i(b_1),
\epsilon_i(b_2)-\alpha_i^\vee(\wt b_1)\}$.
The tensor product of crystals is associative
(cf.~\cite[Proposition~1.3.1]{Ka93}) 
and a tensor product of normal crystals
is also normal (cf. for example~\cite[Lemma~5.2.6]{JB}).

\subl{P35}
Let~$\mathcal{A}$ be the associative monoid generated by the~$e_i,f_i\:
i\in I$. We say that a crystal~$B$ is generated by~$b\in B$ over a
submonoid~$\mathcal A'$ of~$\mathcal{A}$ if~$B=\mathcal A'b:=
\{f b\: f\in\mathcal A'\}\setminus\{0\}$. If~$B$ is generated by~$b$ 
over~$\mathcal{A}$ we will say that~$B$ is generated by~$b$.

Let~$B$ be a crystal. An element~$b\in B$ is said to be of a highest
(respectively, lowest) weight~$\lambda\in P(\pi)$ if~$\wt b=\lambda$ and
$e_i b=0$ (respectively, $f_i b=0$) \forall~$i\in I$. Let~$\mathcal E$
(respectively, $\mathcal F$) be the submonoid of~$\mathcal{A}$ generated
by the~$e_i$ (respectively, by the~$f_i$), $i\in I$. We call~$B$ a
highest (respectively, lowest) weight crystal of highest (respectively,
lowest) weight~$\lambda$ if there exists an element~$b_\lambda$ of
highest (respectively, lowest) weight~$\lambda$ such 
that~$B=\mathcal Fb_\lambda$ (respectively, $B=\mathcal Eb_\lambda$).
\begin{Lemma}
Let~$B$ be a normal crystal and assume that there exists~$b_0\in B$ 
such that~$b_0^{\tensor m}\tensor b\in\mathcal Fb_0^{\tensor m+1}$,
\forall~$b\in B$ and \forall~$m\ge0$. 
Then~$B^{\tensor m}$ is generated by~$b_0^{\tensor m}$ over~$\mathcal F$
\forall~$m>0$.
Similarly, if there exists~$b_0\in B$ such that~$b\tensor b_0^{\tensor m}
\in\mathcal Eb_0^{\tensor m}$ \forall~$b\in B$ 
and~\forall~$m\ge0$ then~$B^{\tensor m}$
is generated by~$b_0^{\tensor m}$ over~$\mathcal E$.
\end{Lemma}
\begin{Proof}
The proof is by induction on~$m$. The induction base is given by the
assumption. Suppose that~$m>1$ and assume that some~$b'\in B^{\tensor m-1}$
satisfies~$b'\tensor b''\in\mathcal
Fb_0^{\tensor m}$ \forall~$b''\in B$. We claim that~$f_ib'\tensor b''
\in\mathcal F
b_0^{\tensor m}$ \forall~$b''\in B$ and \forall~$i\in I$ such that
$f_i b'\not=0$. 
Indeed, $\phi_i(b')>0$ since~$f_i b'\in B^{\tensor m-1}$ and~$B$ is normal.
If~$b''\in B$ satisfies
$\epsilon_i(b'')<\phi_i(b')$, then~$f_ib'\tensor b''=f_i(b'\tensor b'')\in
\mathcal Fb_0^{\tensor m}$ by~\eqref{P30.10b}. Otherwise
set~$k=\epsilon_i(b'')-\phi_i(b')+1$. Then
$0<k\le \epsilon_i(b'')$, whence~$b''':=e_i^k b''\in B$ by normality of~$B$. 
It follows that $b'\tensor b'''\in\mathcal Fb_0^{\tensor m}$. 
On the other hand~$\epsilon_i(b''')=\epsilon_i(b'')-k=
\phi_i(b')-1<\phi_i(b')$, whence
$f_i^{k+1}(b'\tensor b''')=
f_i^k(f_ib'\tensor b''')=f_ib'\tensor f_i^k b'''=
f_ib'\tensor b''$ by~\eqref{P30.10b}. 
Therefore, $f_ib'\tensor b''\in\mathcal Fb_0^{\tensor m}$.

Furthermore, $b_0^{\tensor m-1}\tensor b
\in\mathcal Fb_0^{\tensor m}$ \forall~$b\in B$ by assumption. Then it follows
from the claim by induction on~$k$ that~$f_{i_1}\cdots f_{i_k}b_0^{\tensor 
m-1}\tensor b\in\mathcal Fb_0^{\tensor m}$ \forall~$b\in B$ provided 
that~$f_{i_1}\cdots f_{i_k}b_0^{\tensor m-1}\not=0$. The assertion
follows since~$B^{\tensor m-1}=\mathcal Fb_0^{\tensor m-1}$ by
the induction hypothesis.

Similarly, for the second part it is enough to prove that, for~$b'\in
B^{\tensor m-1}$ fixed, $b''\tensor b'\in\mathcal Eb_0^{\tensor m}$
\forall~$b''\in B$ implies that~$b''\tensor e_ib'\in\mathcal Eb_0^{\tensor m}$
\forall~$b''\in B$ and \forall~$i\in I$ such that~$e_i b'\not=0$ (or,
equivalently, $\epsilon_i(b')>0$).
For, observe that, 
for~$b''\in B$ such that~$\phi_i(b'') < \epsilon_i(b')$, 
$b''\tensor e_ib'=e_i(b''\tensor b')\in\mathcal Eb_0^{\tensor m}$
by~\eqref{P30.10a}. Furthermore, assume that~$\phi_i(b'')\ge \epsilon_i(b')$
and set~$n=\phi_i(b'')-\epsilon_i(b')+1$. Then~$0<n\le
\phi_i(b'')$, whence~$b'''=f_{i}^n b''\in B$ and~$\phi_i(b''')=
\phi_i(b'')-n=\epsilon_i(b')-1<\epsilon_i(b')$.
Then, by~\eqref{P30.10a},
$e_{i}^{n+1}(b'''\tensor b')=e_{i}^n(b'''\tensor e_ib')=
e_{i}^n b'''\tensor e_ib'=b''\tensor e_ib'$, and so
$b''\tensor e_ib'\in\mathcal Eb_0^{\tensor m}$.
\end{Proof}

\subl{P40}
A morphism of crystals~$\psi$ (cf.~\cite[1.2.1]{Ka93}) is a map~$\psi:B_1
\sqcup \{0\}\to B_2\sqcup \{0\}$ such that~$\psi(0)=0$ and,
\forall~$i\in I$,
\begin{enumspecm}{{\mathrm M}_}
\item If~$b\in B_1$ and~$\psi(b)\in B_2$ 
then $\epsilon_i(\psi(b))=\epsilon_i(b)$, $\wt\psi(b)=\wt b$.
\item For all~$b\in B_1$, $\psi(e_ib)=e_i\psi(b)$ provided that~$\psi(e_ib),
\psi(b)\in B_2$.
\item For all~$b\in B_1$, $\psi(f_ib)=f_i\psi(b)$, provided 
that~$\psi(f_ib),\psi(b)\in B_2$.
\end{enumspecm}
A morphism is said to be strict if it commutes with the $e_i,f_i\:i\in I$.
Evidently, any morphism of normal crystals is strict 
([\citen{Ka93}, Lemma~1.2.3]). Throughtout the rest of the paper,
all morphisms of crystals will be assumed to be strict.

Let~$B'$ and~$B$ be crystals. We say that~$B'$ is a subcrystal of~$B$
if there is an injective morphism of crystals from~$B'$ to~$B$. In particular,
a subset $B'\subset B$ will be called a subcrystal of~$B$ 
if $B'$ is a crystal with respect
to the operations $e_i,f_i,\epsilon_i,\phi_i:i\in I$ and~$\wt$ of~$B$
restricted to~$B'$. 
A crystal is said to be indecomposable if it does not admit a
non-empty subcrystal different from itself.

Observe that a crystal~$B$ is indecomposable \iff it is generated
by an element~$b\in B$. Indeed, if~$B$ is indecomposable then 
for any~$b\in B$, $B'=\mathcal{A}b\ni b$ is a non-empty subcrystal of~$B$,
hence coincides with~$B$. On the other hand, suppose that~$B=\mathcal{A}b$
for some~$b\in B$ and that~$B'\subsetneq B$ is a non-empty subcrystal of~$B$.
Then $b\notin B'$ for otherwise~$B=\mathcal Ab\subset B'$. 
On the other hand, for any~$b'\in B'$, there
exists a monomial~$f\in\mathcal A$ 
such that~$b'=fb$. It follows from~\sref{P10}{C}{inv} that there exists a 
monomial~$f'\in\mathcal A$ such that~$b=f'b'$. Since~$B'$ is a crystal, 
we conclude that~$b\in B'$, which is a contradiction. In particular,
it follows that if~$B$ admits a decomposition as a disjoint union of
finitely many indecomposable crystals, then such a decomposition is
unique up to a permutation of the components.

\sect{The crystal~$\hatBlm{}$ and its combinatorics}

\subl{C10}
Set~$I=\{0,\dots,\ell\}$.  Henceforth we identify~$I$ 
with~$\Z/(\ell+1)\Z$ in the sense that $i+k$, $i\in I$, $k\in\Z$ 
is understood as~$i+k\pmod{\ell+1}$.
Let~$A=(a_{ij})_{i,j\in I}$ be the Cartan matrix of the affine Lie algebra
$\lie g=\widehat{\lie{sl}_{\ell+1}}$. 
Explicitly, $a_{ij}=2\delta_{i,j}-\delta_{i,j-1}-\delta_{i-1,j}$,
$i,j\in I$ (for example, $A=\left(\smallmatrix 2&-2\\-2&2\endsmallmatrix
\right)$ for~$\ell=1$).
We will use two different realisations of~$A$.
The first one is the realisation in the sense of~\cite[Chap.~1]{Kac92},
that is, we consider a triple~$(\lie h,\pi,\pi^\vee)$ 
where~$\dim_\Q\lie h=\ell+2$. Observe that~$\delta=\alpha_0+\cdots+
\alpha_\ell\in\lie h^*$ satisfies~$\alpha_i^\vee(\delta)=0$
\forall~$i\in I$. Fix~$\Lambda_i\in\lie h^*$ as in~\ref{P10}. Then
$\Lambda_0,\dots,\Lambda_\ell,\delta$ form a basis of~$\lie h$. 
Throughout the rest of the paper we
take~$P(\pi)=P_0(\pi)\oplus\Z\delta$. The corresponding crystals
will be called affine.

The other realisation is obtained by replacing~$\lie h$ 
by~$\lie h'=\Q\pi^\vee$. 
Then~$\lie h'{}^*=\Q P_0(\pi)\cong \lie h^*/\Q\delta$.
We will use the same notations for
the elements of~$\pi$ and~$\pi^\vee$ in both realisations. The corresponding
crystals will be referred to as finite. The image of the weight map for
finite crystals is contained in~$P_0(\pi)\cong P(\pi)/\Z\delta$.

\subl{C20}
The finite crystal~$B_\ell$ is a set indexed by~$I$. The elements of~$B_\ell$
will be denoted by~$b_i\: i\in I$.
The crystal operators~$e_i,f_i,\epsilon_i,\wt$ on~$B_\ell$ are defined
by the following formulae
\begin{gather}
\begin{alignat*}{3}
&e_i b_j=\delta_{i,j} b_{j-1},&\qquad
&f_i b_j=\delta_{i-1,j} b_i,\\
&\epsilon_i(b_j)=\delta_{i,j},
&&\wt b_j=\Lambda_1-\Lambda_0-\sum_{1\le k<j} 
\alpha_k=-\Lambda_j+\Lambda_{j+1}.
\end{alignat*}
\lbl{10}
\end{gather}
One can easily check that~$B_\ell$ is a normal crystal. Moreover,
if we consider~$B_\ell$ as a crystal with respect to the operations
$e_i,f_i,\epsilon_i\: i\in I\setminus\{0\}$ and the realisation
$(\lie h_0,\pi_0,\pi_0^\vee)$ of~$A_0=(a_{ij})_{i,j\in I\setminus\{0\}}$,
where~$\pi_0=\pi\setminus\{\alpha_0\}$, $\pi_0^\vee
=\pi\setminus\{\alpha_0^\vee\}$ and~$\lie h_0=\Q\pi_0^\vee$,
then~$B_\ell$ can be
realised as a crystal basis of the natural representation~$\C^{\ell+1}$
of the finite dimensional simple Lie algebra~$\lie{sl}_{\ell+1}$
(cf.~\cite{KaN94,Nak98}).

For any~$m\in\N^+$, consider the crystal~$\Blm:=B_\ell^{\tensor m}$,
the operations being defined as in~\ref{P30}. 
Then~$\Blm$ is normal as a tensor product of normal crystals. 
\begin{Lemma}
The crystal~$\Blm:m>0$ is neither a highest weight
nor a lowest weight crystal.
\end{Lemma}
\begin{Proof}
Suppose that there exists~$b\in \Blm$
such that~$e_i b=0$ \forall~$i\in I$. Then, by normality, $\epsilon_i(b)=0$
\forall~$i\in I$. Yet $b$ can be written, uniquely, as~$b=b_j\tensor b'$ 
for some~$j\in I$ and~$b'\in B_\ell^{\tensor m-1}$. Then, 
by~\sref{P30}{T}{eps},
$\epsilon_j(b)=\max\{1,\epsilon_j(b')+1\}>0$, which is a contradiction.
It follows that~$\Blm$ is not a lowest weight crystal either.
Indeed, suppose that~$f_i b=0$ \forall~$i\in I$. Then,
by normality, $\phi_i(b)=0$ \forall~$i\in I$. We claim that
$\sum_{i\in I} \alpha_i^\vee(\wt b)=0$ \forall~$b\in\Blm$. Indeed,
\forall~$j\in I$ one has
$\sum_{i\in I}\alpha_i^\vee(\wt b_j)=
\sum_{i\in I}(\delta_{i,j+1}-\delta_{i,j})=\sum_{i\in I}(\delta_{i-1,j}-
\delta_{i,j})=0$. The claim now follows by~\sref{P30}{T}{wt}.
Then $\sum_{i\in I} \phi_i(b)=\sum_{i\in I} \epsilon_i(b)=0$
by~\sref{P10}{C}{epphi}. Yet~$\epsilon_i(b)\ge0$ by normality of~$\Blm$
and so~$\epsilon_i(b)=0$ \forall~$i\in I$, which is
a contradiction by the first part. 
\end{Proof}

\subl{C30}
Even though~$\Blm$ is not a highest weight crystal, it turns out
to be generated by its element over~$\mathcal F$ or~$\mathcal E$.
\begin{Prop}
The crystal~$\Blm$ is generated by~$b_0^{\tensor m}$ 
over~$\mathcal F$.
Furthermore, $\Blm$ is also generated by~$b_0^{\tensor m}$
over~$\mathcal E$.
\end{Prop}
\begin{Proof}
By~Lemma~\ref{P35} it suffices to prove that~$b_0^{\tensor m-1}\tensor b_i
\in\mathcal Fb_0^{\tensor m}$, \forall~$i\in I$, $m>0$. The cases
$m=1$ and~$m >1$, $i=0$ are trivial. 
Suppose further that~$m>1$ and~$i\not=0$. Since
$\phi_i(b_0^{\tensor m-1})=0=\epsilon_i(b_{i-1})$, 
$i\in I\setminus\{0,1\}$,
$f_i(b_0^{\tensor m-1}\tensor b_{i-1})=b_0^{\tensor m-1}\tensor b_i$
by~\eqref{P30.10b},
whence
$b_0^{\tensor m-1}\tensor b_i=f_i\cdots f_2(b_0^{\tensor m-1}\tensor b_1)$,
$i\in I\setminus\{0,1\}$. Thus, it suffices to prove that~$b_0^{\tensor m-1}
\tensor b_1\in \mathcal Fb_0^{\tensor m}$.
Indeed, observe that~$\phi_i^{}(b_{i-1}^{\tensor k})=k$, whence
$f_i^k b_{i-1}^{\tensor k}=b_i^{\tensor k}$ \forall~$i\in I$, $k>0$.
In particular, $f_1^m b_0^{\tensor m}=b_1^{\tensor m}$. Furthermore,
for all~$i\in I\setminus\{1\}$,
$\phi_i^{}(b_{i-1}^{\tensor m-1})=m-1>\epsilon_i(b_1)=0$, 
whence by~\eqref{P30.10b}~$f_i^{m-1}(b_{i-1}^{\tensor m-1}\tensor b_1)=
f_i^{m-1}b_{i-1}^{\tensor m-1}\tensor b_1=
b_i^{\tensor m-1}\tensor b_1$.
It follows that $b_0^{\tensor m-1}
\tensor b_1 = f_0^{m-1}f_\ell^{m-1}\cdots f_2^{m-1} f_1^m b_0^{\tensor m}
\in \mathcal Fb_0^{\tensor m}$ 
as required.

For the second part, it is sufficient to prove, by Lemma~\ref{P35}, that
$b_i\tensor b_0^{\tensor m-1}\in \mathcal E b_0^{\tensor m-1}$
\forall~$i\in I$ and~$m\ge 1$. Suppose that~$m>1$ and~$i\not=0$,
the other cases being obvious. Since~$e_i(b_i\tensor b_0^{\tensor m-1})=
b_{i-1}\tensor b_0^{\tensor m-1}$, it is sufficient to
prove that~$b_\ell\tensor b_0^{\tensor m-1}\in\mathcal Eb_0^{\tensor m}$
\forall~$m> 1$. 
Indeed, observe that, \forall~$i\in I$ and~$k>0$, 
$\epsilon_i^{}(b_i^{\tensor k})=k$,
whence~$e_i^k b_i^{\tensor k}=b_{i-1}^{\tensor k}$.
In particular, $e_0^m b_0^{\tensor m}=b_\ell^{\tensor m}$. 
Furthermore, \forall~$i\in I\setminus\{0\}$,
$\epsilon_i^{}(b_i^{\tensor m-1})=m-1>
\phi_i(b_\ell)=0$, whence by~\eqref{P30.10a}~$e_i^{m-1}(b_\ell^{}
\tensor b_i^{\tensor m-1})=b_\ell^{}\tensor e_i^{m-1} b_i^{\tensor 
m-1}=
b_\ell^{}\tensor b_{i-1}^{\tensor m-1}$.
It follows that~$b_\ell\tensor b_0^{\tensor m-1}=
e_1^{m-1}\cdots e_\ell^{m-1} e_0^m b_0^{\tensor m}\in
\mathcal Eb_0^{\tensor m}$ as required.
\end{Proof}
\begin{Rem}
One can show that~$B_\ell$ is a perfect crystal (cf.~\cite[4.6]{KKMMNN92}). 
Then $B_\ell^{\tensor m}$ is indecomposable for all~$m>0$
by~\cite[Corollary~4.6.3]{KKMMNN92}. However, we need a stronger
version of this result, namely that~$B_\ell^{\tensor m}$ is generated
by some element over~$\mathcal F$ and not just over~$\mathcal A$. 
Besides, our proof does not use the fact that~$B_\ell$ is perfect.
\end{Rem}

\subl{C35}
The affine crystal~$\hatBlm{}$, which we are about to define, 
provides the affinisation of~$\Blm$ in the sense of~\cite[3.3]{KKMMNN92}. 
Set
$\hatBlm{}=\Blm\times\Z$ and define 
the crystal operations  as follows.
Denote an element~$(b,n)\:
b\in\Blm$, $n\in\Z$ by~$b\tensor z^n$.
Then~$\wt(b\tensor z^n)=\wt b+n\delta\in P(\pi)$, which is compatible
with the decomposition~$P(\pi)=P_0(\pi)\oplus\Z\delta$.
Furthermore, set
\begin{align*}
&e_i(b\tensor z^n)=
\begin{cases}e_ib\tensor z^{n+\delta_{i,0}},&\text{if~$e_ib\in\Blm$}\\
0,&\text{if~$e_ib=0$}\end{cases}
\\
&f_i(b\tensor z^n)=
\begin{cases}f_ib\tensor z^{n-\delta_{i,0}},&\text{if~$f_ib\in\Blm$} \\
0,&\text{if~$f_ib=0$}\end{cases}
\end{align*}
and~$\epsilon_i(b\tensor
z^n)=\epsilon_i(b)$, $i\in I$. 
Evidently, $\hatBlm{}$ is a normal crystal.  
\begin{Prop}
The crystal~$\hatBlm{}$ is the disjoint union of indecomposable
normal subcrystals~$\hatBlm n$, $n=0,\dots,m-1$, where~$\hatBlm k\:k\in\Z$
is the subcrystal of~$\hatBlm{}$ generated by~$b_0^{\tensor m}\tensor z^k$.
\end{Prop}
The sections~\ref{C37}--\ref{C50} are devoted to the prove of the
above Proposition.

\subl{C37}
The~$\hatBlm n$ are indecomposable by~\ref{P40}
and normal as subcrystals of a normal crystal. So, it 
remains to prove that~$\hatBlm r=\hatBlm s$ if~$s=r\pmod m$,
$\hatBlm r\intersection\hatBlm s=\emptyset$ otherwise, and that
every element of~$\hatBlm{}$ lies in some~$\hatBlm k$. 
\begin{Lemma}
The crystal~$\hatBlm{}$ is a union of~$\hatBlm n\:n=0,\dots,m-1$.
\end{Lemma}
\begin{Proof}
Let~$b\in \Blm$. By Proposition~\ref{C30}, there exists~$f=
f_{i_1}\cdots f_{i_k}\in\mathcal F$ such 
that~$b=fb_0^{\tensor m}$. Define
$$
n_f(b):=\# \{ t \: i_t = 0\}.
$$
Then~$b\tensor z^r=f(b_0^{\tensor m}\tensor z^{r+n_f(b)})\in\hatBlm{r+
n_f(b)}$.

Furthermore, an elementary computation shows that
\begin{align*}
&(e_1^m \cdots e_\ell^m e_0^m)^r (b_0^{\tensor m}
\tensor z^k)=b_0^{\tensor m}\tensor z^{k+rm},\\
&(f_0^m f_\ell^m\cdots f_1^m)^r(b_0^{\tensor m}\tensor z^k)=
b_0^{\tensor m}\tensor z^{k-rm},\qquad r>0.
\end{align*}
It follows immediately that~$\hatBlm r=\hatBlm s$ if~$r=s\pmod m$.
\end{Proof}

\subl{C40}
Notice that~$n_f(b)$ depends on~$f$. For example, one has
$f_0^m f_\ell^m\cdots f_1^m b_0^{\tensor m}=b_0^{\tensor m}$.
However, it turns out that the residue class of~$n_f(b)$ modulo~$m$ 
does not depend on~$f$, which allows one to introduce a function
$N:\hatBlm{}\to\Z/m\Z$ such that~$N(b)=n$ \iff~$b\in\hatBlm{n}$. That
function also plays a crucial role in the computation of characters
of the indecomposable subcrystals of~$\hatBlm{}$ and in the construction
of a subcrystal of Littelmann's path crystal isomorphic to~$\hatBlm{}$.

Given a product~$u=b_{j_k}\tensor \cdots\tensor b_{j_1}$, 
define $t(u)=j_k$, $h(u)=j_1$ and~$|u|=k$. 
\begin{defi}
Let $b=b_{i_m}\tensor\cdots\tensor b_{i_1}$ be an element of~$\Blm$.
Define 
\begin{equation}\lbl{10}
\maj(b):=\{ r\: 1\le r<m,\, i_r > i_{r+1}\}. 
\end{equation}
Furthermore, set~$k=\#\maj(b)+1$ and write~$\maj(b)=\{
n_1<\cdots < n_{k-1}\}$ if~$k>1$. Set~$n_0=0$, $n_k=|b|=m$ and define
\begin{align*}
\tmaj (b)&:=\{n_0,\dots,n_k\}\\
N(b)&:=\sum_{r=1}^{k} r(n_r-n_{r-1}).
\end{align*}
\end{defi}
\begin{Prop}
Let~$b\in \Blm$. Then
\begin{enumeraa}
\item $N(e_ib)=N(b)-\delta_{i,0}\pmod m$ provided that~$e_ib\in\Blm$.
\item $N(f_ib)=N(b)+\delta_{i,0}\pmod m$ provided that~$f_ib\in\Blm$.
\end{enumeraa}
\end{Prop}
\begin{Proof*}
Observe that the second statement follows from the first.
Indeed, if $b'=f_ib\in B$, then~$b=e_ib'$ by~\sref{P10}{C}{inv}, 
whence~$N(f_ib)=N(b')=N(e_ib')+\delta_{i,0} \pmod m=N(b)+\delta_{i,0}\pmod m$.

Suppose that~$b=b''\tensor b_i\tensor b'$ for
some~$b'$, $b''$, possibly empty. We claim that if~$e_i b=b''
\tensor b_{i-1}\tensor b'$ then~$h(b'')\not=i-1$ and~$t(b')\not=i$.
Indeed, suppose that~$b''=b'''\tensor b_{i-1}$ and~$e_i b=
b'''\tensor b_{i-1}^{\tensor 2}\tensor b'$. Then by~\sref{P30}{T}{e} we must
have, in particular,~$\epsilon_i(b_{i-1})-\alpha_i^\vee(\wt b''')<
\epsilon_i(b_i)-\alpha_i^\vee(\wt b''')-\alpha_i^\vee(\wt b_{i-1})$. That
inequality reduces to~$0<1-\alpha_i^\vee(\wt b_{i-1})=1-
\alpha_i^\vee(\Lambda_{i}-\Lambda_{i-1})=0$, 
which is a contradiction. Similarly, if~$t(b')=i$, that is~$b'=
b_i\tensor b'''$, and~$e_i b=b''\tensor b_{i-1}\tensor
b_i\tensor b'''$ then we have, 
by~\sref{P30}{T}{e}, $\epsilon_i(b_i)-\alpha_i^\vee(\wt b'')\ge
\epsilon_i(b_i)-\alpha_i^\vee(\wt b'')-\alpha_i^\vee(\wt b_i)=
\epsilon_i(b_i)-\alpha_i^\vee(\wt b'')+1$, which is absurd.

Suppose that~$e_i b\not=0$. Since~$e_i b_j=0$ if~$j\not=i$,
$b=b''\tensor b_i\tensor b'$ for some~$b'$, $b''$ such that 
$e_i b=b''\tensor e_i b_i\tensor b'=b''\tensor b_{i-1}\tensor b'$.
First, consider the case~$i\not=0$. Since~$h(b'')\not=i-1$
and~$t(b')\not=i$ by the above, 
we conclude that~$\maj(b)=\maj(e_i b)$ whence~$N(e_ib)=N(b)$.

Suppose now that~$i=0$ and retain the notations from the above definition.
\begin{enumera}
\item Assume first that~$b=b_0\tensor b'$ and~$e_0 b=b_\ell\tensor b'$. 
By the above claim, $t(b')\not=0$. 
Then~$n_{k-1}=m-1$ and~$\tmaj (e_0 b)=\{n_0,\dots,n_{k-2},n_k\}$.
Therefore,
\begin{equation*}
\begin{split}
N(e_0 b)&=\sum_{r=1}^{k-2} r(n_r-n_{r-1})+(k-1)(n_k-n_{k-2})=
\sum_{r=1}^{k-1} r(n_r-n_{r-1})+k-1\\&=
\sum_{r=1}^{k-1} r(n_r-n_{r-1})+k(n_k-n_{k-1})-1
=N(b)-1.
\end{split}
\end{equation*}
\item Suppose that~$b=b''\tensor b_0\tensor b'$, where~$|b'|, |b''|>0$,
and~$e_0 b=b''\tensor b_\ell\tensor b'$.
Since~$t(b')\not=0$ by the above claim, $|b'|\in\maj(b)$. 
Suppose that~$n_s=|b'|$ in our notations for
the elements of~$\tmaj(b)$. On the other hand,
$h(b'')<\ell$ by the claim we proved above. Since~$\ell\ge t(b')$,
it follows that
$\tmaj (e_0 b)=\{n_0,\dots,n_{s-1},n_s+1,n_{s+1},\dots, n_k\}$, 
whence
\begin{equation*}
\begin{split}
N(e_0 b)&=\sum_{r=1}^{s-1} r(n_r-n_{r-1})+s(n_s-n_{s-1}+1)+
(s+1)(n_{s+1}-n_s-1)+\sum_{r=s+2}^k r(n_r-n_{r-1})\\
&=\sum_{r=1}^k r(n_r-n_{r-1})-1=N(b)-1.
\end{split}
\end{equation*}
\item Finally, assume that~$b=b''\tensor b_0$ and~$e_0 b=b''\tensor b_\ell$.
Evidently, $1\notin\maj(b)$. On the other hand,
$h(b'')<\ell$ by our claim, whence~$\tmaj(e_0b)=\{n_0,1,n_1,\dots,n_k\}$.
Therefore,
\begin{align*}
N(e_0 b)&=1+2(n_1-1)+\sum_{r=2}^k (r+1)(n_r-n_{r-1})=
-1+\sum_{r=1}^k (r+1)(n_r-n_{r-1})\\&=-1+N(b)+\sum_{r=1}^k (n_r-n_{r-1})
=N(b)+m-1=N(b)-1\pmod m.\tag*{\qed}
\end{align*}
\end{enumera}
\end{Proof*}
\begin{cor}
Let~$b$ be an element of~$\Blm$ and let~$f\in \mathcal F$
be a monomial such that~$b=fb_0^{\tensor m}$. Then~$n_f(b)=N(b)\pmod m$.
In particular, the residue class of~$n_f(b)$ modulo~$m$ does not depend
on~$f$.
\end{cor}
\begin{Proof}
It suffices to observe that~$N(b_0^{\tensor m})=m=0\pmod m$.
\end{Proof}

\subl{C50}
Now we are able to complete the proof of Proposition~\ref{C35}.
\begin{Proof}
By Lemma~\ref{C37}, $\hatBlm{}$ is a union of~$\hatBlm n\:n=0,\dots,m-1$.
It only remains to prove that~$\hatBlm{r}\intersection\hatBlm{s}$ is empty
if~$r\not=s\pmod m$.
Given~$b=b'\tensor z^k$, $b'\in\Blm$, $k\in\Z$, 
set~$N(b)=N(b')+k$ and define
$$
C_n:=\{b\in\hatBlm{}\: N(b)=n\pmod m\}.
$$
Evidently, $C_r\intersection C_s$ is empty if~$r\not=s\pmod m$.
The idea is to prove that~$C_n=\hatBlm{n}$. 

First, let us prove
that~$C_n$ is a subcrystal of~$\hatBlm{}$. Indeed,
let~$b=b'\tensor z^k\in C_n$ and suppose that~$e_ib\not=0$. 
Then~$e_ib=e_ib'\tensor
z^{k+\delta_{i,0}}$ and~$e_ib'\not=0$. It follows from Proposition~\ref{C40}
that $N(e_ib')=N(b')-\delta_{i,0}\pmod m$. Therefore,
$N(e_ib)=N(e_ib')+k+\delta_{i,0}=N(b)\pmod m$, hence~$e_ib\in
C_n$. Similarly, if~$f_ib\not=0$, then~$f_ib'\not=0$ 
and~$N(f_ib')=N(b')+\delta_{i,0}\pmod m$, 
whence~$N(f_ib)=N(f_ib')+k-\delta_{i,0}=N(b)\pmod m$.

Furthermore, $N(b_0^{\tensor m})=0
\pmod m$ hence~$C_n$ contains~$b_0^{\tensor m}\tensor z^n$. 
By the proof of
Lemma~\ref{C37}, the $b_0^{\tensor m}\tensor z^{n+rm}$ 
lie in the subcrystal of~$\hatBlm{}$ 
generated by~$b_0^{\tensor m}\tensor z^n$ \forall~$r\in\Z$ which is
contained in~$C_n$. Take $b=b'\tensor
z^k\in C_n$ and let~$f\in\mathcal F$ be a monomial such that~$b'=fb_0^{
\tensor m}$. Then~$b\tensor z^k=f(b_0^{\tensor m}\tensor z^{k+n_f(b')})$.
On the other hand, $k+n_f(b')=k+N(b')\pmod m=N(b)\pmod m=n$ by 
Corollary~\ref{C40}. Thus, $C_n$ is generated by~$b_0^{\tensor m}\tensor z^n$,
hence is indecomposable by~\ref{P40}, 
and contains~$\hatBlm{n}$ as a subcrystal by the definition of the latter.
\end{Proof}
\begin{cor}
The indecomposable subcrystals of~$\hatBlm{}$ are given explicitly as
$$
\hatBlm{n}=\{b\tensor z^k\: b\in\Blm,\, k\in\Z,\, N(b)=n-k\pmod m\},
\quad n=0,\dots,m-1.
$$
\end{cor}
\begin{Rem}
The decomposition of~$\hatBlm{}$ of Proposition~\ref{C35} 
appears in~\cite[Corollary~6.25]{Nak98}
in the case~$\ell=1$. However,
our proof for arbitrary~$\ell$ does not use the theory of perfect crystals, 
yields an efficient explicit description of the indecomposable subcrystals 
and allows one to compute their formal characters.
\end{Rem}

\subl{C70}
Our present aim is to compute the formal character of~$\hatBlm n$.
For, we calculate first the cardinalities of the
sets~$\Blm_\nu^n:=\{ b\in \Blm\: \wt b=\nu, N(b)=n\pmod m\}$,
where~$\nu\in P_0(\pi)$ and~$n=0,\dots,m-1$. Let~$b=b_{i_m}\tensor
\cdots\tensor b_{i_1}$ be an element of~$\Blm$ and set
$k_i=\#\{r\: i_r=i\}$. Then~$\wt b=\sum_{i\in I} k_i\wt b_i=
\sum_{i\in I} k_i(\Lambda_{i+1}-\Lambda_i)$.
On the other hand, the numbers
$k_i\: i\in I$ are uniquely determined by~$\wt b$ and~$m$. Indeed, write
$\wt b=\sum_{i\in I} k_i \wt b_i$. Then~$\alpha_i^\vee(\wt b)=
k_{i-1}-k_i$, $i\in I\setminus\{0\}$, whence~$k_i=
k_\ell+\sum_{j>i} \alpha_j^{\vee}(\wt b)$, $i\in I\setminus\{\ell\}$. 
Yet $\sum_{i\in I} k_i=m$, whence~$(\ell+1)k_\ell=
m-(\alpha_1^\vee(\wt b)+2\alpha_2^\vee(\wt b)+\cdots+\ell\alpha_\ell^\vee(
\wt b))$. 

Thus, there is a bijection between the set of weights of~$\Blm$
and the set~$\{(k_0,\dots,k_\ell)\in\N^{\ell+1}\: 
\sum_{i\in I} k_i = m\}$. We will
identify a weight~$\nu$ of~$\Blm$ with the tuple~$(k_0,\dots,k_\ell)$.
It follows immediately that 
$\# \Blm_\nu=\binom{m}{k_0,\dots, k_\ell}$. Indeed, the
multinomial coefficient~$\binom{m}{k_0,\dots,k_\ell}$ 
gives the number of distinct permutations
the word~$0^{k_0}\cdots\ell^{k_\ell}$. 
By the above there is a bijection between this set and the set of
all elements of~$\Blm$ of weight~$\nu=(k_0,\dots,k_\ell)$.
\begin{Prop}
Let~$\nu=(k_0,\dots,k_\ell)$ be a weight of~$\Blm$. Then
$$
\# \Blm_\nu^n:=\frac1m \sum_{d|\gcd(k_0,\dots,k_\ell)}
\phi_{-n}(d)\binom{\frac md}{\frac{k_0}d,\dots, \frac{k_\ell}d},\qquad
\phi_r(d)=\phi(d)\,\frac{\mu(d/\gcd(d,r))}{\phi(d/\gcd(d,r))},
$$
where~$\phi$ is the Euler function and~$\mu$ is the M\"obius function,
$\mu(k)=0$ if~$k$ is divisible by a square and~$\mu(k)=(-1)^r$ if~$k$
is a product of~$r$ distinct primes.
\end{Prop}
\noindent
The proof of this proposition is not based on the theory of
crystals, and for that reason is given in the Appendix.

\subl{C80}
Retain the notations of~\ref{I20}. 
\begin{Thm}
The formal character of~$\hatBlm{n}$, $n=0,\dots,m-1$ 
equals that of the simple
integrable module~$L(\ell,m;n)$ described in~\ref{I20}.
\end{Thm}
\begin{Proof}
Evidently, $\Omega(\hatBlm n)\subset \{\nu+k\delta\:
\nu\in\Omega(\Blm),\,k\in\Z\}$. 
Let~$\nu=\sum_{i=0}^\ell k_i(\Lambda_{i+1}-\Lambda_i)$ be a weight of~$\Blm$.
By Corollary~\ref{C50}, $\hatBlm{n}_{\nu+k\delta}=
\{b\tensor z^k\: b\in \Blm_\nu,\, N(b)=n-k\pmod m\}$, whence
$$
\# \hatBlm{n}_{\nu+k\delta}=\#\Blm_\nu^{n-k} =
\frac 1m\,\sum_{d|\gcd(k_0,\dots,k_\ell)} 
\phi_{k-n}(d)\binom{\frac md}{\frac{k_0}d,\dots,
\frac{k_\ell}d}
$$
by Proposition~\ref{C70}.
On the other hand, $\nu+k\delta$ is a weight of~$L(\ell,m;n)$
and the dimension of the corresponding weight space equals the 
right hand side of the above expression by~\cite[Theorem~4.4]{Grxx}.
\end{Proof}

\sect{Littelmann's path crystal and~$\hatBlm{}$}

\subl{L10}
Let us briefly recall the definition of Littelmann's path 
crystal~([\citen{Li94,Li95}]). 

Fix a realisation~$(\lie h,\pi,\pi^\vee)$ of a symmetrizable
Cartan matrix~$A$ and denote by~$[a,b]$ the set~$\{\tau\in\Q\:a\le\tau\le b\}$.
Let~$\mathbb P$ be the set of piecewise-linear continuous paths~$b:[0,1]\to
\lie h$ such that~$b(0)=0$ and~$b(1)\in P(\pi)$. Two paths~$b_1$,~$b_2$ are
considered to be identical if there exists a piecewise-linear,
nondecreasing, surjective continuous map~$\phi:[0,1]\to[0,1]$ such
that~$b_1=b_2\circ\phi$. 

One can endow~$\mathbb P$ with a structure of a normal crystal in the following
way. For all~$i\in I$, define the Littelmann function
$h^i_b:[0,1]\to\Q$, $h^i_b(\tau)=-\alpha^\vee_i(b(\tau))$ and 
set~$\epsilon_i(b)=\max\{h^i_b(\tau)\intersection\Z\:\tau\in[0,1]\}$. 
Furthermore,
define~$e^i_+(b)=\min\{\tau\in[0,1]\: h^i_b(\tau)=\epsilon_i(b)\}$. 
If~$e^i_+(b)=0$, define~$e_i b =0$. Otherwise 
let~$e^i_-(b)=\max\{\tau\in[0,e^i_+(b)]\:
h^i_b(\tau)=\epsilon_i(b)-1\}$ and define
$$
(e_ib)(\tau):=\begin{cases}
b(\tau),&\tau\in[0,e^i_-(b)]\\
s_i(b(\tau)-b(e^i_-(b)))+b(e^i_-(b)),&\tau\in[e^i_-(b),e^i_+(b)]\\
b(\tau)+\alpha_i,&\tau\in[e^i_+(b),1],
\end{cases}
$$
where~$s_i$ is the simple reflection corresponding to~$\alpha_i$,
$s_i\lambda=\lambda-\alpha_i^\vee(\lambda)\alpha_i^{}$ \forall~$\lambda
\in\lie h^*$ and~$s_i(b(\tau))$ is taken point-wise. In particular,
$s_i(b(\tau))=b(\tau)+h^i_b(\tau)\alpha_i$.

Similarly, define~$f^i_+(b)=\max\{\tau\in[0,1]\: h^i_b(\tau)=\epsilon_i(b)\}$.
If~$f^i_+(b)=1$, set~$f_ib=0$. Otherwise, let~$f^i_-(b)=
\min\{\tau\in[f^i_+(b),1]\: h^i_b(\tau)=\epsilon_i(b)-1\}$ and define
$$
(f_ib)(\tau):=\begin{cases}
b(\tau),&\tau\in[0,f^i_+(b)]\\
s_i(b(\tau)-b(f^i_+(b)))+b(f^i_+(b)),&\tau\in[f^i_+(b),f^i_-(b)]\\
b(\tau)-\alpha_i,&\tau\in[f^i_-(b),1].
\end{cases}
$$
Finally, set~$\wt b(\tau)=b(1)$.
\begin{Rem}
We use the definition of crystal operations on~$\mathbb P$ 
given in~\cite[6.4.4]{JB} which differs by the sign of~$h^i_b$ from
the original definition of~\cite[1.2]{Li94}.
That choice is more convenient for the proof of Proposition~\ref{L20}.
\end{Rem} 
A path~$b\in\mathbb P$ is said to have the {\em integrality property}
(cf.~\cite[2.6]{Li95}) if the maximal value of~$h^i_b(\tau)$
is an integer \forall~$i\in I$. If that condition holds for every $b\in B$
where~$B$ is a subcrystal of~$\mathbb P$, we say that~$B$ has the
integrality property.

\subl{L15}
For any~$b_1,b_2\in\mathbb P$, let~$b_1 * b_2$ denote their
concatenation, that is, a path defined by
$$
(b_1 * b_2)(\tau)=\begin{cases} b_1(\tau/\sigma),&\tau\in[0,\sigma]\\
b_1(1)+b_2((\tau-\sigma)/(1-\sigma)),&\tau\in[\sigma,1],
\end{cases}
$$
where~$\sigma\in(0,1)$. One may check that the resulting path does not
depend on~$\sigma$, up to a reparametrisation. Moreover,
the concatenation of paths is compatible with the tensor product rules 
listed in~\ref{P30} (cf.~\cite[2.6]{Li95}).
Henceforth we will use the notation~$b_1\tensor b_2$
for the concatenation of~$b_1$, $b_2$.
 
\subl{L20}
Retain the notations of~\ref{C10}--\ref{C50}.
Our present aim is to define an isomorphism between~$\hatBlm{}$ and
a certain subcrystal of Littelmann's path crystal~$\mathbb P$. 

Let~$\blambda=(\lambda_0,\dots,\lambda_r)$ be a tuple of
elements of~$\lie h^*=\Q P(\pi)$. 
We assume that~$\lambda_0=0$ and~$\lambda_r\in
P(\pi)$. Given $\mathbf a=(a_0,\dots,a_r)$, $a_s\in\Q$ such that
$0=a_0<a_1<\cdots < a_r=1$, define a path~$p_{\blambda,
\mathbf a}(\tau)\in\mathbb P$ as
$$
p_{\blambda,\mathbf a}(\tau)=\lambda_{r-1}+\left(\frac{\tau-a_{r-1}}{a_r-
a_{r-1}}\right)(\lambda_r-\lambda_{r-1}),\qquad
\tau\in[a_{r-1},a_r].
$$
Evidently~$p_{\blambda,
\mathbf a}\in\mathbb P$. We shall omit~$\mathbf a$ if~$a_s=s/r$, 
$s=0,\dots,r$. 

Let~$b=b_{i_m}\tensor\cdots\tensor b_{i_1}$ be an element of~$\Blm$.
We associate to~$b\tensor z^n$ a path~$p_{\blambda}$, 
$\blambda=(\lambda_0,\dots,\lambda_m)$, where~$\lambda_s=
\sum_{t=m-s+1}^m \wt b_{i_t}+ \kappa_s(b\tensor z^n)\delta$, $s=0,\dots,m$ and 
the~$\kappa_s(b\tensor z^n)$
are defined recursively in the following way. Set~$\kappa_0(b\tensor z^n)=0$.
Furthermore, write~$\tmaj(b)=\{n_0,\dots,n_k\}$ as in Definition~\ref{C40}
and let~$\rho_s(b)$ be the unique~$r$, $1\le r\le k$ such that
$n_{r-1}< m-s+1\le n_r$.
Then
\begin{equation}\lbl{5}
\kappa_s(b\tensor z^n) = \kappa_{s-1}(b\tensor z^n)-(\rho_s(b)-1)+(N(b)+n-m)/m, \qquad s=1,\dots,m.
\end{equation}

\begin{Prop}
The map~$\psi: \hatBlm{}\to\mathbb P$ given by~$b\tensor z^n\mapsto
p_{\blambda}$ with~$\blambda$ defined as above is an
injective morphism of normal crystals and the image of~$\psi$
has the integrality property.
\end{Prop}
\begin{Proof}
The injectivity of~$\psi$ is obvious. 
Let~$b=b_{i_m}\tensor\cdots\tensor b_{i_1}\in
\Blm$. By~\sref{P30}{T}{eps},
\begin{equation*}
\begin{split}
\epsilon_i(b)&=\max_{1\le s\le m}\{\epsilon_i(b_{i_s})-\sum_{t> s}
\alpha_i^\vee(\wt b_{i_t})\}=
\max_{1\le s\le m}
\{\delta_{i,i_s}-\alpha_i^\vee(\lambda_{m-s}-\kappa_{m-s}(b\tensor z^n)\delta)\}\\
&=
\max_{1\le s\le m}\{\delta_{i,i_s}-\alpha_i^\vee(\lambda_{m-s})\}.
\end{split}
\end{equation*}
On the other hand,
\begin{equation*}
\begin{split}
h^i_{p_{\blambda}}(\tau)&=-\alpha_i^\vee(\lambda_{s-1})-
(m\tau-s+1)\alpha_i^\vee(\wt b_{i_{m-s+1}})\\&=
-\alpha_i^\vee(\lambda_{s-1})+(m\tau-s+1)(\delta_{i,i_{m-s+1}}-
\delta_{i-1,i_{m-s+1}}),\qquad\tau\in[(s-1)/m,s/m].
\end{split}
\end{equation*}
It follows immediately that~$\max\{h^i_{p_\blambda}(\tau)\:
\tau\in[(s-1)/m,s/m]\}=\delta_{i,i_{m-s+1}}-\alpha_i^\vee(\lambda_{s-1})\in\Z$.
Since~$h^i_{p_{\blambda}}$ is linear on the intervals~$[(t-1)/m,t/m]$,
$t=1,\dots,m$ and all the local maxima are integral, we conclude
that~$p_\blambda$ has the integrality property and that the maximal (integer) 
value of~$h^i_{p_{\blambda}}$ is attained at~$t/m$ for some~$t$. 
Therefore,
\begin{equation}\lbl{10}
\epsilon_i(p_{\blambda})=\max_{1\le s\le m}\{\delta_{i,i_{m-s+1}}-
\alpha_i^\vee(\lambda_{s-1})\}=\max_{1\le s\le m}\{\delta_{i,i_s}-
\alpha_i^\vee(\lambda_{m-s})\}=
\epsilon_i(b)=\epsilon_i(b\tensor z^n).
\end{equation}
Thus, $\epsilon_i$ commutes with~$\psi$ \forall~$i\in I$. 
Furthermore, by the definition of~$\kappa_s(b\tensor z^n)$ 
$$
\kappa_m(b\tensor z^n)=-\sum_{r=1}^k (r-1)(n_r-n_{r-1})+N(b)-m+n=n,
$$
whence~$\wt p_{\blambda}=p_{\blambda}(1)=
\sum_{s=1}^m \wt b_{i_s}+n\delta=\wt b+n\delta=\wt(b\tensor z^n)$.
Thus, $\psi$ commutes with~$\wt$ and hence with the~$\phi_i$ 
\forall~$i\in I$. 

Since both~$\hatBlm{}$ and~$\mathbb P$ are normal, 
it remains to prove that~$\psi$ commutes with the~$e_i$, $f_i\:i\in I$.
By~\sref{P10}{C}{inv} it suffices to prove that~$\psi$ commutes with the~$e_i$.
Suppose that~$e_i b\not=0$. Then~$e_ib=b_{i_m}\tensor\cdots
\tensor e_i b_{i_s}\tensor \cdots\tensor b_{i_1}$ and~$i_s=i$. 
By~\sref{P30}{T}{e} $s$ is
the largest element of the set~$\{1,\dots,m\}$ such that
$\epsilon_i(b_{i_s})-\sum_{t>s}\alpha_i^\vee(\wt b_{i_t})=\epsilon_i(b)$.
It follows from~\eqref{L20.10} and~\ref{L10} that
$e^i_+(p_{\blambda})=(m-s+1)/m$ and~$e^i_-(p_{\blambda})=(m-s)/m$. 
Using the definition of~$e_i$ given in~\ref{L10}, we obtain
$$
e_i p_{\blambda}=p_{\blambda'},
$$
where
\begin{equation}\lbl{15}
\lambda'_t=\begin{cases}\lambda_t,&
t=0,\dots,m-s,\\
\lambda_t+\alpha_i,& t=m-s+1,\dots,m.
\end{cases}
\end{equation}
Let us prove that~$\blambda'=\bmu$ 
where~$p_\bmu=\psi(e_i b\tensor z^{n+\delta_{i,0}})$. Indeed,
since~$e_ib=b_{i_m}\tensor \cdots\tensor b_{i_{s+1}}
\tensor b_{i-1}\tensor b_{i_s-1}\tensor\cdots
 b_{i_1}$
and~$\wt b_{i-1}-
\wt b_i=2\Lambda_i-\Lambda_{i-1}-\Lambda_{i+1}$ which equals
$\alpha_i-\delta_{i,0}\delta$ as an element of~$P(\pi)$, 
it follows from the definition
of the~$\lambda_t$ and~\eqref{L20.15} that
\begin{equation}\lbl{15b}
\mu_t=\begin{cases}\lambda'_t+(k'_t-k_t)\delta,& t=0,\dots,m-s,\\
\lambda'_t+(k'_t-k_t-\delta_{i,0})\delta,& t=m-s+1,\dots,m,
\end{cases}
\end{equation}
where $k'_t=\kappa_t(e_i(b\tensor z^n))=
\kappa_t(e_ib\tensor z^{n+\delta_{i,0}})$ and~$k_t=\kappa_t(b\tensor z^n)$.
Denote also~$r_t=\rho_t(b)$, $r'_t=\rho_t(e_i b)$.

Suppose first that~$i\not=0$. Then~$N(e_ib)=N(b)$ and, as we saw in the
proof of~Proposition~\ref{C40}, $\maj(e_ib)=\maj(b)$, whence~$r'_t=r_t$
\forall~$t=1,\dots,m$. 
It follows that~$k'_t=k_t$ \forall~$t=0,\dots,m$ and so~$\bmu=\blambda'$.

Assume that~$i=0$. First, if~$s=1$, that is $e_0 b=b_{i_m}
\tensor\cdots\tensor e_0 b_{i_1}$, then
by the proof of Proposition~\ref{C40},
$N(e_0b)=N(b)+m-1$ and~$\maj(e_0b)=\maj(b)\sqcup\{1\}$.
It follows that~$r'_t=r_t+1$, $t=1,\dots,m-1$. Then
\begin{equation*}
\begin{split}
k'_t&=k'_{t-1}-(r'_t-1)+(N(e_0b)-m+n+1)/m=k'_{t-1}-r_t+(N(b)+n)/m
\\
&=k'_{t-1}-(r_t-1)+(N(b)-m+n)/m,
\end{split}
\end{equation*}
whence~$k'_t=k_t$ \forall~$t=0,\dots,m-1$.
Furthermore, $r'_m=1=r_m$, whence
$$
k'_m=k'_{m-1}+(N(b)+n)/m=k_{m-1}+(N(b)+n)/m=k_m+1.
$$
Then~$\bmu=\blambda'$ by~\eqref{L20.15b}.

Finally, assume that~$s>1$. Then, by the proof of Proposition~\ref{C40},
$N(e_0b)=N(b)-1$ and
$r'_t=r_t$, $t\not=m-s+1$ whilst $r'_{m-s+1}=r_{m-s+1}-1$. One has
\begin{equation}\lbl{20}
\begin{split}
k'_t&=k'_{t-1}-(r'_t-1)+(N(e_0b)+n+1-m)/m\\
&=k'_{t-1}-(r_t-1)+(N(b)+n-m)/m,
\end{split}
\end{equation}
\forall~$t=1,\dots,m-s,m-s+2,\dots,m$. It follows that
$k'_t=k_t$ \forall~$t=0,\dots,m-s$. Furthermore,
\begin{equation*}
\begin{split}
k'_{m-s+1}&=k'_{m-s}-(r'_{m-s+1}-1)+(N(b)+n-m)/m\\&=k_{m-s}-(r_{m-s+1}-1)+1
+(N(b)+n-m)/m
=k_{m-s+1}+1.
\end{split}
\end{equation*}
Then one concludes, using the recurrence~\eqref{L20.20} for the~$k'_t$ and 
the definition of the~$k_t$ that~$k'_t=k_t+1$, $t=m-s+1,\dots,m$. 
Therefore, $\bmu=\blambda'$ by~\eqref{L20.15b}.
\end{Proof}

\subl{L25}
Our main result (Theorem~\ref{I20}) follows immediately from
Theorem~\ref{C80} and
\begin{Thm}
The indecomposable crystal~$\hatBlm n$ is isomorphic to a subcrystal
of the Littelmann's crystal~$\mathbb P$ generated by the path
$p_{\ell,m,n}(\tau):=(m(\Lambda_1-\Lambda_0)+n\delta)\tau$, $\tau\in[0,1]$.
\end{Thm}
\begin{Proof}
Let~$p_\blambda$ be the image of~$b_0^{\tensor m}\tensor z^n$ under
the map constructed in~\ref{L20}. Then~$\blambda=
(\lambda_0,\dots,\lambda_m)$, where~$\lambda_t=t(\Lambda_1-\Lambda_0)+(nt/m)
\delta$, $t=0,\dots,m$. Since~$(\lambda_t-\lambda_{t-1})/(a_t-a_{t-1})=
m(\Lambda_1-\Lambda_0+(n/m)\delta)=m(\Lambda_1-\Lambda_0)+n\delta$, it
follows that~$p_\blambda$ coincides with~$p_{\ell,m,n}$. On the
other hand, $b_0^{\tensor m}\tensor z^n$ generates an indecomposable
subcrystal~$\hatBlm{n}$ of~$\hatBlm{}$, whence~$p_\blambda$ generates
an indecomposable subcrystal of~$\mathbb P$, which is isomorphic 
to~$\hatBlm{n}$
by Proposition~\ref{L20}.
\end{Proof}

\sect{Decomposition of~$B(\Lambda)\tensor\hatBlm{}$}

\subl{L30}
Let~$P^+(\pi)=\{\lambda\in P(\pi)\: \alpha_i^\vee(\lambda)\ge0,\,\forall
i\in I\}$ and~$\mathbb P^+=\{b\in\mathbb P\: \alpha_i^\vee(b(\tau))\ge0,
\,\forall i\in I,\,\forall\tau\in[0,1]\}$. Take~$\lambda\in P^+(\pi)$.
A path~$b\in\mathbb P$ is
said to be $\lambda$-dominant (cf.~\cite{Li94})
if~$\alpha_i^\vee(\lambda+b(\tau))\ge0$
\forall~$\tau\in[0,1]$ and \forall~$i\in I$.
The following Lemma is
rather standard (cf. for example~\cite[6.4.14]{JB})
\begin{Lemma}
Let~$\lambda\in P^+(\pi)$ and $b_\lambda\in \mathbb P^+$ such that
$\wt b_\lambda=\lambda$. Let~$b\in\mathbb P$ and
suppose that~$b$ has the integrality property. Then the following
are equivalent
\begin{enumer}
\item $\epsilon_i(b_\lambda\tensor b)=0$ \forall~$i\in I$.
\item $\epsilon_i(b)\le \alpha_i^\vee(\lambda)$ \forall~$i\in I$.
\item $b$ is $\lambda$-dominant.
\item $b_\lambda\tensor b\in\mathbb P^+$.
\end{enumer}
\end{Lemma}
\begin{Proof}
Suppose that~{(i)} holds. 
By~\sref{P30}{T}{eps}, $0=\epsilon_i(b_\lambda\tensor b)=
\max\{\epsilon_i(b_\lambda),\epsilon_i(b)-\alpha_i^\vee(\lambda)\}=
\max\{0,\epsilon_i(b)-\alpha_i^\vee(\lambda)\}$
implies that~$\epsilon_i(b)\le\alpha_i^\vee(\lambda)$, hence~{(ii)}
follows from~{(i)}. Suppose further that~$\epsilon_i(b)\le \alpha_i^\vee(b)$.
Recall that~$\epsilon_i(b)=\max\{-\alpha_i^\vee(b(\tau))\intersection \Z\:
\tau\in[0,1]\}=\max\{-\alpha_i^\vee(b(\tau))\:\tau\in[0,1]\}$ since~$b$
is assumed to have the integrality property. Therefore, {(ii)} 
implies~{(iii)}. Furthermore, $b_\lambda\tensor b=b_\lambda(\tau/\sigma)$,
$\tau\in[0,\sigma]$, whence~$\alpha_i^\vee((b_\lambda\tensor b)(\tau))\ge0$
\forall~$\tau\in[0,\sigma]$. On the other hand, $b_\lambda\tensor b=
b_\lambda(1)+b((\tau-\sigma)/(1-\sigma))$, $\tau\in[\sigma,1]$,
whence~$\alpha_i^\vee((b_\lambda\tensor b)(\tau))=
\alpha_i^\vee(\lambda)+\alpha_i^\vee(b(\tau'))\ge0$, where
$\tau'=(\tau-\sigma)/(1-\sigma)\in[0,1]$, provided that $b$ is
$\lambda$-dominant. Thus, {(iii)} leads to~{(iv)}. Finally, if~$b_\lambda
\tensor b\in\mathbb P^+$, then $-\alpha_i^\vee((b_\lambda\tensor b)(\tau))
\le 0$ \forall~$i\in I$ and \forall~$\tau\in[0,1]$, whence the maximum of 
$h^i_{b_\lambda\tensor b}$ is non-positive. Since~$h^i_{b_\lambda\tensor b}(0)
=0$, we conclude that~$\epsilon_i(b_\lambda\tensor b)=0$. Thus, {(iv)}
implies~{(i)}.
\end{Proof}

\subl{L35}
Take~$\lambda\in P^+(\pi)$ and let~$b_\lambda\in\mathbb P^+$ be any path 
satisfying~$\wt b_\lambda=\lambda$.
By the Isomorphism Theorem
of Littelmann (cf.~\cite[Theorem~7.1]{Li95}),
the subcrystal of~$\mathbb P$ generated by~$b_\lambda$ is isomorphic
to the subcrystal~$B(\lambda)$ of~$\mathbb P$ generated by the
path~$\tau\mapsto \lambda\tau$.
Moreover, by~\cite[7, Corollary 1~b)]{Li95}, $b_\lambda$ is the unique highest
weight element of the subcrystal of~$\mathbb P$ it generates. Observe
that if~$\lambda\in\Z\delta$, then the corresponding crystal~$B(\lambda)$
is trivial.

Let~$B$ be a subcrystal of~$\mathbb P$. Given~$\lambda\in
P^+(\pi)\setminus\Z\delta$, let $B^\lambda$ be the set 
of~$\lambda$-dominant paths in~$B$.
\begin{Lemma}
Let~$B$ be a subcrystal of~$\mathbb P$ and suppose that~$B$ has 
the integrality property. Then
$$
B(\lambda)\tensor B\longtwoheadrightarrow \coprod_{b\in B^\lambda}
B(\lambda+\wt b)
$$
and the crystals which appear in the right-hand side are the only
highest weight subcrystals of~$B(\lambda)\tensor B$.
\end{Lemma}
\begin{Proof}
Let~$b_\lambda=\lambda\tau$, $\tau\in[0,1]$. 
By Lemma~\ref{L30}{(i)}, $b_\lambda\tensor b$, where~$b\in B^\lambda$,
is a highest weight element and, by the Isomorphism 
Theorem~\cite[Theorem~7.1]{Li95} and Lemma~\ref{L30}{(iv)} it generates
a highest weight subcrystal
of~$\mathbb P$ isomorphic to~$B(\mu)$ where~$\mu=(b_\lambda\tensor b)(1)=
\lambda+\wt b$. It remains to prove that
if $b'\tensor b\in B(\lambda)\tensor B$ lies entirely in the dominant
Weyl chamber then~$b'=b_\lambda$ and~$b\in B^\lambda$. For,
assume that~$b'\tensor b\in\mathbb P^+$. Then~$\epsilon_i(b'\tensor
b)=0$ by Lemma~\ref{L30}. On the other hand, $\epsilon_i(b'\tensor b)
=\max\{\epsilon_i(b'),\epsilon_i(b)-\alpha_i^\vee(\wt b')\}\ge
\epsilon_i(b')$ \forall~$i\in I$ by~\sref{P30}{T}{eps}.
Therefore, $\epsilon_i(b')=0$ \forall~$i\in I$, whence
$b'=b_\lambda$ by~\cite[7, Corollary~1~b)]{Li95}. It remains
to apply Lemma~\ref{L30}{(iii)}.
\end{Proof}

\subl{L40}
Henceforth, let~$B$ be the image of~$\hatBlm{}$ inside~$\mathbb P$
under the morphism~$\psi$ constructed in~\ref{L20}. By~Proposition~\ref{L20}
$B$ has the integrality property, hence we immediately obtain 
a surjective morphism of~$B(\lambda)\tensor B$ onto a disjoint
union of highest weight crystals~$B(\mu)$ where~$\mu=\lambda+\wt b$
for some~$b\in B^\lambda$. 
Our goal now is to prove that this surjective morphism 
is actually an isomorphism. By Lemma~\ref{L35}, 
it suffices to prove that~$B(\lambda)
\tensor B$ is generated by its highest weight elements over~$\mathcal F$
and that~$B^\lambda$ is not empty for all~$\lambda\in P^+(\pi)
\setminus\Z\delta$.
\begin{Lemma} Let~$b$ be an element of~$\Blm$. Then
$$
\epsilon_j(b\tensor b_i)=\begin{cases}
\epsilon_j(b)+1,&\text{if~$i=j$ and~$\phi_i(b)=0$}\\
\epsilon_j(b),&\text{otherwise.}
\end{cases}
$$
In particular, if~$b\in\Blm$ then there exists~$i\in I$ such that
$\epsilon_j(b\tensor b_i)=\epsilon_j(b)$ \forall~$j\in I$.
\end{Lemma}
\begin{Proof}
By~\sref{P10}{C}{epphi} and~\sref{P30}{T}{eps}, 
$\epsilon_j(b\tensor b_i)=\max\{\epsilon_j(b),
\epsilon_j(b_i)-\alpha_j^\vee(\wt b)\}=\epsilon_j(b)+
\max\{0,\delta_{i,j}-\phi_j(b)\}$. 
The first statement follows
immediately since~$\Blm$ is a normal crystal. The second statement
follows from the first and Lemma~\ref{C20}.
\end{Proof}

\subl{L45}
Set~$\Blm^\lambda=\{b\in\Blm\: \epsilon_i(b)\le\alpha_i^\vee(\lambda),\,
\forall i\in I\}$.
The next step is to prove that every element of~$\Blm$ can
be transformed into an element of~$\Blm^\lambda$ provided
that the latter is not empty by applying some
special monomial~$e\in\mathcal E$.
\begin{Lemma}
Let~$\lambda\in P(\pi)^+\setminus\Z\delta$ and $b\in\Blm$.
Then there exist $i_1,\dots,i_k\in I$ such that
$e_{i_k}^{m_k}\cdots e_{i_1}^{m_1} b\in\Blm^\lambda$ and
$m_r=\epsilon_{i_r}(u_{r-1})-\alpha_{i_r}^\vee(\lambda)>0$,
where~$u_0=b$ and~$u_r=e_{i_r}^{m_r} u_{r-1}^{}$, $r=1,\dots,k$.
\end{Lemma}
\begin{Proof}
Write~$b=b_{j_1}\tensor\cdots\tensor b_{j_m}$ and choose~$1\le s\le m$
maximal such that~$b'=b_{j_1}\tensor\cdots\tensor b_{j_s}$ satisfies
$\epsilon_j(b')\le\alpha_j^\vee(\lambda)$ \forall~$j\in I$. If~$s=m$
then there is nothing to prove. Otherwise, write~$b=b'\tensor b_i\tensor b''$
where~$i=j_{s+1}$.
By Lemma~\ref{L40}, $\epsilon_j(b'\tensor b_i)=
\epsilon_j(b')\le\alpha_j^\vee(\lambda)$ \forall~$j\not=i$. On the
other hand, $\epsilon_j(b'\tensor b_i)>\alpha_j^\vee(\lambda)$
for some~$j\in I$ by the choice of~$b'$. It follows from Lemma~\ref{L40}
that~$j=i$, $\epsilon_i(b'\tensor b_i)=\epsilon_i(b')+1$,
$\epsilon_i(b')=\alpha_i^\vee(\lambda)$ and~$\phi_i(b')=0$. In particular,
$e_i(b'\tensor b_i)=b'\tensor b_{i-1}$ by~\eqref{P30.10a}.
Furthermore,
\begin{equation*}
\begin{split}
\epsilon_i(b)&=\max\{\epsilon_i(b'\tensor b_i),
\epsilon_i(b'')-\alpha_i^\vee(\wt b')-\alpha_i^\vee(\wt b_i)\}
\\&=\epsilon_i(b')+1+\max\{0,\epsilon_i(b'')-\phi_i(b')\}
=\alpha_i^\vee(\lambda)+1+\max\{0,\epsilon_i(b'')\}.
\end{split}
\end{equation*}
Set~$i_1=i$. Then~$m_1=\epsilon_i(b)-\alpha_i^\vee(\lambda)=
\max\{0,\epsilon_i(b'')\}+1>0$. Furthermore, by~\sref{P30}{T}{e},
$$
e_{i_1}^{m_1}b=e_i(b'\tensor b_i)\tensor e_i^{\epsilon_i(b'')}b''=
b'\tensor b_{i-1}\tensor u_1''.
$$
If~$\epsilon_j(b'\tensor b_{i-1})=\epsilon_j(b')\le\alpha_j^\vee
(\lambda)$ for all~$j\in I$
then we can use induction on~$|b'|$. Otherwise, we repeat
the above argument for~$r=1,\dots,k$ where~$k\in\N$ is minimal
such that~$i-k$ satisfies~$\epsilon_j(b'\tensor b_{i-k})=
\epsilon_j(b')\le\alpha_j^\vee(\lambda)$ \forall~$j\in I$. The
existence of such~$i-k\in I$ is guaranteed by Lemma~\ref{L40}.
As a result we obtain a monomial~$e=e_{i_k}^{m_k}\cdots e_{i_1}^{m_1}$
where~$m_{i_r}=\epsilon_{i_r}(u_{r-1})-\alpha_{i_r}^\vee(\lambda)$,
$u_r=e_{i_r}^{m_r}\cdots e_{i_1}^{m_1}b$ and~$u_0=b$ such that~$e b 
= \tilde b'\tensor u_k''$ where~$\epsilon_j(\tilde b')
\le \alpha_j^\vee(\lambda)$ \forall~$j\in I$ and~$|\tilde b'|>
|b'|$. The assertion follows by induction on~$|b'|$.
\end{Proof}

\subl{L50}
The following Lemma allows one to prove that~$B^\lambda$ is not empty
for each~$\lambda\in P^+(\pi)\setminus\Z\delta$.
\begin{Lemma}
For all $i\in I$, $\Blm^{\Lambda_i}=\{b(i,m)\}$ where~$b(i,m)=b_i\tensor
b_{i+1}\tensor\cdots\tensor b_{i+m-1}$. 
\end{Lemma}
\begin{Proof}
One has, \forall~$t=1,\dots,m$,
$$
r^j_t(b(i,m))=\epsilon_j(b_{i+t-1})-\sum_{s=1}^{t-1}
\alpha_j^\vee(\wt b_{i+s-1})=\delta_{j,i+t-1}+\delta_{j,i}-\delta_{j,i+t-1}
=\delta_{i,j},
$$
whence~$\epsilon_j(b(i,m))=\max_t\{r^j_t(b(i,m))\}=\delta_{i,j}$.
Therefore, $b(i,m)\in \Blm^{\Lambda_i}$.

We prove that~$b(i,m)$ is the only element of~$\Blm^{\Lambda_i}$
by induction on~$m$. The induction
base is given by~\eqref{C20.10}. 
Suppose that~$b\in\Blm\: m >1$ satisfies~$\epsilon_j(b)
\le\delta_{i,j}$. Then~$\epsilon_j(b)=\delta_{i,j}$ for otherwise
$b$ is a highest weight element of~$\Blm$ by normality of the latter,
which is a contradiction by Lemma~\ref{C20}.
Write $b=b_r\tensor 
b'$ where $b'\in B_\ell(m-1)$. Since~$\epsilon_j(b)=\max\{\epsilon_j(b_r),
\epsilon_j(b')+\delta_{j,r}-\delta_{j,r+1}\}=
\max\{0,\epsilon_j(b')-\delta_{j,r+1}\}+\delta_{j,r}$
it follows that~$\epsilon_r(b) > 0$, whence~$r=i$.
Then $\epsilon_j(b)=\delta_{j,i}+\max\{0,\epsilon_j(b')-\delta_{j,i+1}\}$,
hence we must also have~$\epsilon_j(b')=0$ if~$j\not=i+1$ 
and~$\epsilon_{i+1}(b')
\le 1$. Therefore, $b'\in B_\ell(m-1)^{\Lambda_{i+1}}$, whence
$b'=b(i+1,m-1)$ by the induction hypothesis, and so~$b=b_i\tensor b(i+1,
m-1)=b(i,m)$.
\end{Proof}

\subl{L55}
Now we are able to prove Theorem~\ref{I30}
\begin{Proof}
By Proposition~\ref{L20} and Lemma~\ref{L30}, $B^\lambda=\{ \psi(b\tensor 
z^k)\: b\in\Blm^\lambda,\, k\in\Z\}$. Yet \forall~$\lambda\in P^+(\pi)
\setminus\Z\delta$,
there exists $i\in I$ such that~$\alpha_i^\vee(\lambda)>0$. It follows
that~$\Blm^{\Lambda_i}$, which is non-empty by Lemma~\ref{L50}, is
contained in~$\Blm^\lambda$. Therefore, $B^\lambda$ is not empty for 
all~$\lambda\in P^+(\pi)\setminus\Z\delta$.

By~\ref{L40} and Lemma~\ref{L35}, it remains to prove that~$B(\lambda)
\tensor B$ is generated over~$\mathcal F$ by its highest weight 
elements $b_\lambda\tensor p$ where~$p\in B$ is $\lambda$-dominant
That is equivalent to proving that \forall~$b\in B(\lambda)$ and 
\forall~$p\in B$
there exists~$e\in\mathcal E$ such that~$e(b\tensor p)=
b_{\lambda}\tensor p'$ where~$p'$ is~$\lambda$-dominant.
By~Lemma~\ref{L30} and Proposition~\ref{L20}, 
$p'$ is $\lambda$-dominant \iff~$p'=\psi(b'\tensor z^k)$
where~$b'\in\Blm^\lambda$ and~$k\in\Z$.

Take arbitrary~$b\in B(\lambda)$, $p\in B$ and let
us prove first that there exists a monomial~$e\in\mathcal E$
such that~$e(b\tensor p)=b_{\lambda}\tensor p'$
for some~$p'\in B$. Indeed, by~\eqref{P30.10a},
$e_j^{k+1}(b\tensor p)=e_jb\tensor e_j^kp=e_jb\tensor p'$, 
where $k=\max\{0,\epsilon_j(p)-\phi_j(b)\}\le \epsilon_j(p)$. 
Since~$B(\lambda)$ is generated
by~$b_{\lambda}$ over~$\mathcal F$ by~\cite[7, Corollary~1~c)]{Li95}, 
there exists a monomial~$e\in\mathcal E$
such that~$eb=b_{\lambda}$. By the above, there exists a monomial~$e'$
such that~$e'(b\tensor p)=eb\tensor p'=b_{\lambda}\tensor p'$ 
and~$p'\not=0$.

It remains to prove that for all~$p\in B$, there exist~$e\in
\mathcal E$ such that~$e(b_{\lambda}\tensor p)=
b_{\lambda}\tensor p'$ where~$p'\in B^\lambda$.
Set~$m_j(p)=\max\{0,\epsilon_j(p)-\alpha_j^\vee(\lambda)\}$. Then
\begin{equation}\lbl{10}
e_j^{m_j(p)}(b_{\lambda}\tensor p)=b_{\lambda}\tensor e_j^{m_j(p)}p.
\end{equation}
Indeed, $\epsilon_j(b_{\lambda}\tensor p)=\max\{0,\epsilon_j(p)-
\alpha_j^\vee(\lambda)\}=m_j(p)$ and,
since~$e_jb_{\lambda}=0$, we conclude that~$e_j^{m_j(p)}(b_{\lambda}\tensor
p)=b_{\lambda}\tensor e_j^{m_j(p)}p$. 
Furthermore,
write~$p=\psi(u\tensor z^n)$.
By Lemma~\ref{L45}, there exists a monomial~$e=e_{i_k}^{m_k}
\cdots e_{i_1}^{m_1}$ such that~$e(u\tensor z^n)=u'\tensor z^{n+s}$
where~$u'\in\Blm^\lambda$,
$s=\sum_t m_t \delta_{0,i_t}$ and~$m_t=m_{i_t}(u_{t-1})>0$ where~$u_0=u$,
$u_t=e_{i_t}^{m_t}\cdots e_{i_1}^{m_1}u$, $t=1,\dots,k$.
Then, using~\eqref{L55.10} repeatedly, we
obtain 
\begin{equation*}
e(b_{\lambda}\tensor p)=e_{i_k}^{m_k}\cdots
e_{i_2}^{m_2}(b_{\lambda}\tensor e_{i_1}^{m_1}p)=
\cdots=b_{\lambda}\tensor ep=b_{\lambda}\tensor p',
\end{equation*}
where~$p'=\psi(u'\tensor z^{n+s})\in B^\lambda$.
\end{Proof}

\subl{L60}
In the case~$\lambda=\Lambda_i$ we are able to describe the decomposition 
of Theorem~\ref{I30} more explicitly.
\begin{Prop}
For all~$i\in I$, 
$$
B(\Lambda_i)\tensor \hatBlm{}\stackrel{\sim}{\to} \coprod_{k\in\Z}
B(\Lambda_{i+m}+k\delta).
$$
\end{Prop}
\begin{Proof}
Observe that~$\wt b(i,m)=\Lambda_{i+m}-\Lambda_i$. The assertion
follows immediately from Theorem~\ref{I30} and Lemma~\ref{L50}.
\end{Proof}

\appendix*{Proof of Proposition~\ref{C70}}

\def\Maj{\operatorname{Maj}}

\subl{A10}
Given~$b\in\Blm$, define its major index~$\Maj(b):=\sum_{r=1}^{k-1} n_r$
where~$\tmaj(b)=\{n_0<\cdots<n_k\}$ (cf.~Definition~\ref{C40}). 
In other words, $\Maj(b)$ equals 
the sum of all elements in~$\maj(b)$ and zero if~$\maj(b)$ is empty.
Our definition of~$\Maj(b)$ is just the standard definition of
the major index of a word in a free monoid over a completely ordered alphabet.
\begin{Lemma}
For all~$b\in\Blm$, $\Maj(b)=-N(b)\pmod m$.
\end{Lemma}
\begin{Proof*}
Indeed, write~$\tmaj(b)=\{n_0,\dots,n_k\}$ and recall that~$n_0=0$, $n_k=m$. 
Then
\begin{align*}
N(b)&=\sum_{r=1}^k r(n_r-n_{r-1})=
\sum_{r=1}^k r n_r - \sum_{r=1}^{k-1} (r+1)n_r\\
&= kn_k-\sum_{r=1}^{k-1} n_r
=km-\Maj(b)=-\Maj(b)\pmod m.
\tag*{\qed}
\end{align*}
\end{Proof*}

\subl{A20}
Set~$[n]=(q^n-1)/(q-1)=1+q+\cdots+q^{n-1}$ and define
$$
[n]!:=[n] [n-1]\cdots [1],
\qquad \qbinom{n}{n_1,\dots, n_k}=\frac{[n]!}{[n_1]!\cdots
[n_k]!},
$$
where~$n=n_1+\cdots+n_k$. It is convenient to assume
that~$\qbinom{n}{n_1,\dots, 
n_k}=0$ if~$n\not=n_1+\cdots+n_k$. We will also use the notation~$[n]_{q^r}
:=(q^{nr}-1)/(q^r-1)=1+q^r+\cdots+q^{(n-1)r}$.
\begin{Lemma}
The cardinality of the set~$\Blm_\nu^{-n}=\{b\in \Blm_\nu\: N(b)=-n\pmod m\}$
equals the coefficient of~$q^n$ in the polynomial
$$
\qbinom{m}{k_0,\dots, k_\ell}\pmod{q^m-1},
$$
where~$\nu=(k_0,\dots,k_\ell)$.
\end{Lemma}
\begin{Proof}
Let~$\Gamma=\{\gamma_1,\dots,\gamma_r\}$
be a completely ordered alphabet and let~$R$ be a set of all distinct
permutations of the word~$\gamma_1^{k_1}\cdots \gamma_r^{k_r}$ in the
free monoid over~$\Gamma$. Then
$$
\sum_{w\in R} q^{\Maj(w)}=
\qbinom{m_1+\cdots+m_r}{m_1,\dots, m_r}
$$
by a classical theorem of MacMahon (cf. for example~\cite[Theorem~3.7]{AnB}). 
Apply this theorem to~$\Gamma=\{b_\ell,\dots,b_0\}$ and
the word~$b_\ell^{k_\ell}\tensor\cdots\tensor b_0^{k_0}$ whose distinct
permutations form the set~$\Blm_\nu$ where~$\nu=\sum_{i\in I} 
k_i (\Lambda_{i+1}-\Lambda_i)$. The result then follows immediately
from Lemma~\ref{A10}.
\end{Proof}

\subl{A30}
Let~$r_0,\dots,r_\ell$ be non-negative integers and denote their sum by~$r$.
Set, for all~$n\in\Z$,
$$
C(r_0,\dots,r_\ell;n):=\frac1r \sum_{d|r} \phi_n(d)\binom{\frac rd}{\frac{r_0}d
,\dots, \frac{r_\ell}d},
$$
where~$\phi_n(d)$ is defined as in Proposition~\ref{C70}. 
Furthermore, for all~$d$ dividing~$r$ set
$$
\tilde C(r_0,\dots,r_\ell;d):=\frac{d}{r}\sum_{dd'|r}
\mu(d')\binom{\frac{r}{dd'}}{\frac{r_0}{dd'},\dots,
\frac{r_\ell}{dd'}}.
$$
From now on we adopt the convention that a multinomial coefficient
equals zero if any of the rational numbers involved is not an integer.
\begin{Lemma}
Fix~$m\in\N^+$ and let~$k_0,\dots,k_\ell$ be non-negative integers
such that~$k_0+\cdots+k_\ell=m$. Then
\begin{equation}\lbl{10}
\sum_{n=0}^{m-1} C(k_0,\dots,k_\ell;n)q^n=
\sum_{d|m}\tilde C(k_0,\dots,k_\ell;d) [m/d]_{q^d}.
\end{equation}
\end{Lemma}
\begin{Proof}
The coefficient of~$q^n$ in the right-hand side of~\eqref{A30.10}
equals
\begin{align*}
\sum_{d|n}\tilde C(k_0,\dots,k_\ell;d)&=
\sum_{d|n}\sum_{dd'|m}\frac dm\,\mu(d')\binom{\frac{m}{dd'}}{\frac{k_0}{dd'}
,\dots,\frac{k_\ell}{dd'}}\\
&=\frac1m\,\sum_{d|m} \binom{\frac{m}{d}}{\frac{k_0}d,\dots,
\frac{k_\ell}d} \Big(
d \sum_{d'|d,\, d|nd'} \frac{\mu(d')}{d'}\Big).
\end{align*}
The inner sum equals~$\phi_n(d)$ (cf. the proof of~\cite[Corollary~4.2]{Gr99}).
\end{Proof}
\subl{A40}
The next step of our proof is the following
\begin{Lemma}
Let~$d$ be a divisor of~$m$ and denote by
$\Phi_d(q)$ the~$d$th cyclotomic polynomial. 
Let~$\psi_d$ be the canonical projection~$\Q[q]\to\Q[q]/(\Phi_d(q))
\cong \Q(\zeta_d)$, where~$\zeta_d$ is a~$d$th primitive root of unity. Then
$$
\sum_{r|m} \tilde C(k_0,\dots,k_\ell;r) \psi_d([m/r]_{q^r})=
\binom{\frac md}{\frac{k_0}d,\dots,\frac{k_\ell}d}.
$$
\end{Lemma}
\begin{Proof}
Set~$P(q):=\sum_{r|m} \tilde C(k_0,\dots,k_\ell;r)[m/r]_{q^r}$.
Then~$\psi_d(P(q))$ 
identifies with~$P(\zeta_d)$. 
If~$d$ divides~$r$, then~$\psi_d([m/r]_{q^r})=
m/r$. Otherwise, $\zeta_d^r\not=1$ and so~$\psi_d([m/r]_{q^r})=
(\zeta_d^m-1)/(\zeta_d^r-1)=0$ since~$d|m$. Therefore,
$$
P(\zeta_d)=\sum_{rd|m} \tilde C(k_0,\dots,k_\ell; rd)\, \frac{m}{rd}=
\sum_{r|\frac md} \tilde C\big(\tfrac{k_0}d,\dots,\tfrac{k_\ell}d; r\big) 
\frac{m}{rd}.
$$
Furthermore, set~$n=m/d$, $n_i=k_i/d$, $i=0,\dots,\ell$. Then
$$
P(\zeta_d)=\sum_{r|n}\sum_{d'r|n}\mu(d')\binom{\frac{n}{d'r}}{%
\frac{n_0}{d'r},\dots,\frac{n_\ell}{d'r}}=
\sum_{r|n} \binom{\frac nr}{\frac{n_0}r,\dots,\frac{n_\ell}r}
\sum_{r'|r}\mu(r/r').
$$
It remains to apply the fundamental property of the M\"obius function,
namely, $\sum_{d|n}\mu(d)=0$ if~$n>1$ and~$1$ otherwise.
\end{Proof}

\subl{A50}
The following Lemma has been adapted from~\cite[Lemma~34.1.2]{LB}.
We deem it necessary to present its proof here since we use 
the definition~\cite[3.3]{AnB} of the~$q$-multinomial coefficients,
which differs from that of~\cite[1.3]{LB} by a power of~$q$.
\begin{Lemma}[{cf.~\cite[Lemma~34.1.2]{LB}}]
Let~$m$, $d$ be non-negative integers and let~$\psi_d$ be the map defined
in~\ref{A40}.
\begin{enumera}
\item Let~$k_1,\dots,k_r\in\N$ and
suppose that~$d$ does not divide~$\gcd(k_1,\dots,k_r)$. Then
$$
\psi_d\Big(\qbinom{md}{k_1,\dots, k_r}\Big)=0.
$$
\item For all~$k_1,\dots,k_r\in\N$
$$
\psi_d\Big(\qbinom{md}{k_1 d,\dots, k_r d}\Big)=
\binom{m}{k_1,\dots, k_r}.
$$
\end{enumera}
\end{Lemma}
\begin{Proof}
First, let us prove both assertions for~$r=2$. Set~$\qbinom{m}{k}:=
\qbinom{m}{k,m-k}$.
\begin{enumera}
\item
Obviously, $\qbinom{md}{k}=0$ if~$m=0,1$ and~$d$
does not divide~$k$. Suppose that the assertion holds 
for all non-negative integers~$<m$ and for all~$k < (m-1)d$ not
divisible by~$d$. Take~$k < md$ not divisible by~$d$.
Then by~\cite[Theorem~3.4]{AnB} or~\cite[1.3.1]{LB},
$$
\psi_d\Big(\qbinom{md}{k}\Big)=\sum_{t=0}^k \psi_d(q^{(md-d-t)(k-t)})
\psi_d\Big(\qbinom{(m-1)d}{t}\Big)\psi_d\Big(\qbinom{d}{k-t}\Big)=0,
$$
since at least one of~$t$, $k-t$ is not divisible by~$d$.
\item Since~$\qbinom{md}{kd}=0$ if~$k>m$, we immediately
conclude that the second assertion holds for~$m=0,1$. Furthermore, 
\begin{align*}
\psi_d\Big(\qbinom{md}{kd}\Big)&=\sum_{t=0}^{kd} \psi_d(q^{(m-d-t)(kd-t)})
\psi_d\Big(\qbinom{(m-1)d}{t}\Big)\psi_d\Big(\qbinom{d}{kd-t}\Big)
\\
&=\sum_{t=0}^k \psi_d(q^{(m-d-td)(k-t)d})
\psi_d\Big(\qbinom{(m-1)d}{td}\Big)\psi_d\Big(\qbinom{d}{(k-t)d}\Big),
\end{align*}
where we applied the first part. Notice that~$\psi_d(q^d)=1$. Then,
by the induction hypothesis,
$$
\psi_d\Big(\qbinom{md}{kd}\Big)=\sum_{t=0}^k \binom{m-1}{t}\binom{1}{k-t}=
\binom{m-1}{k}+\binom{m-1}{k-1}=\binom{m}{k}.
$$
\end{enumera}
Suppose now that~$r>2$ and observe that
$$
\qbinom{m}{k_1,\dots, k_r}=\qbinom{m}{k_1} \qbinom{m-k_1}{k_2,\dots, k_r}.
$$
The assertion follows immediately by induction on~$r$.
\end{Proof}

\subl{A60}
Now we are able to prove Proposition~\ref{C70}.
\begin{Proof}
By Lemma~\ref{A20}, $\#\Blm_\nu^{-n}$ equals the coefficient of~$q^n$
in the polynomial~$\qbinom{m}{k_0,\dots, k_\ell}\pmod{q^m-1}$.
It follows from Lemmata~\ref{A30},~\ref{A40} and~\ref{A50} that
$$
\psi_d\Big(\sum_{n=0}^{m-1} C(k_0,\dots,k_\ell;n)q^n\Big)=
\psi_d\Big(\qbinom{m}{k_0,\dots, k_\ell}\Big),
$$
for all~$d$ dividing~$m$. On the other hand,~$q^m-1=\prod_{d|m}\Phi_d(q)$.
Since cyclotomic polynomials are irreducible over~$\Q$
and~$\Q[q]$ is a unique factorisation
domain, it follows that
$$
\sum_{n=0}^{m-1} C(k_0,\dots,k_\ell;n) q^n=
\qbinom{m}{k_0,\dots, k_\ell}\pmod{q^m-1}.
$$
Therefore, $\#\Blm_\nu^{-n}=C(k_0,\dots,k_\ell;n)$.
\end{Proof}


\bibliographystyle{amsplain}
\ifx\undefined\bysame
\newcommand{\bysame}{\leavevmode\hbox to3em{\hrulefill}\,}
\fi

\end{document}